# CONFIDENCE SETS FOR SPLIT POINTS IN DECISION TREES

By Moulinath Banerjee [1] and Ian W. McKeague[2]

*University of Michigan and Columbia University*

We investigate the problem of finding confidence sets for split points in decision trees (CART). Our main results establish the asymptotic distribution of the least squares estimators and some associated residual sum of squares statistics in a binary decision tree approximation to a smooth regression curve. Cube-root asymptotics with nonnormal limit distributions are involved. We study various confidence sets for the split point, one calibrated using the subsampling bootstrap, and others calibrated using plug-in estimates of some nuisance parameters. The performance of the confidence sets is assessed in a simulation study. A motivation for developing such confidence sets comes from the problem of phosphorus pollution in the Everglades. Ecologists have suggested that split points provide a phosphorus threshold at which biological imbalance occurs, and the lower endpoint of the confidence set may be interpreted as a level that is protective of the ecosystem. This is illustrated using data from a Duke University Wetlands Center phosphorus dosing study in the Everglades.

**1. Introduction.** It has been over twenty years since decision trees (CART) came into widespread use for obtaining simple predictive rules for the classification of complex data. For each predictor variable $X$ in a (binary) regression tree analysis, the predicted response splits according to whether $X \leq d$ or $X > d$, for some split point $d$. Although the rationale behind CART is primarily statistical, the split point can be important in its own right, and in some applications it represents a parameter of real scientific interest. For example, split points have been interpreted as thresholds for the presence of environmental damage in the development of pollution control standards. In

Received December 2003; revised March 2006.
[1]Supported by NSF Grant DMS-03-06235.
[2]Supported by NSF Grant DMS-05-05201.
*AMS 2000 subject classifications.* 62G08, 62G20, 62E20.
*Key words and phrases.* CART, change-point estimation, cube-root asymptotics, empirical processes, logistic regression, Poisson regression, nonparametric regression, split point.

a recent study [16] of the effects of phosphorus pollution in the Everglades, split points are used in a novel way to identify threshold levels of phosphorus concentration that are associated with declines in the abundance of certain species. The present paper introduces and studies various approaches to finding confidence sets for such split points.

The split point represents the best approximation of a binary decision tree (piecewise constant function with a single jump) to the regression curve $E(Y|X = x)$, where $Y$ is the response. Bühlmann and Yu [3] recently studied the asymptotics of split-point estimation in a homoscedastic nonparametric regression framework, and showed that the least squares estimator $\hat{d}_n$ of the split point $d$ converges at a cube-root rate, a result that is important in the context of analyzing bagging. As we are interested in confidence intervals, however, we need the exact form of the limiting distribution, and we are not able to use their result due to an implicit assumption that the "lower" least squares estimator $\hat{\beta}_l$ of the optimal level to the left of the split point converges at $\sqrt{n}$-rate (similarly for the "upper" least squares estimator $\hat{\beta}_u$). Indeed, we find that $\hat{\beta}_l$ and $\hat{\beta}_u$ converge at *cube-root rate*, which naturally affects the asymptotic distribution of $\hat{d}_n$, although not its rate of convergence.

In the present paper we find the joint asymptotic distribution of $(\hat{d}_n, \hat{\beta}_l, \hat{\beta}_u)$ and some related residual sum of squares (RSS) statistics. Homoscedasticity of errors is not required, although we do require some mild conditions on the conditional variance function. In addition, we show that our approach readily applies in the setting of generalized nonparametric regression, including nonlinear logistic and Poisson regression. Our results are used to construct various types of confidence intervals for split points. Plug-in estimates for nuisance parameters in the limiting distribution (which include the derivative of the regression function at the split point) are needed to implement some of the procedures. We also study a type of bootstrap confidence interval, which has the attractive feature that estimation of nuisance parameters is eliminated, albeit at a high computational cost. Efron's bootstrap fails for $\hat{d}_n$ (as pointed out by Bühlmann and Yu [3], page 940), but the subsampling bootstrap of Politis and Romano [14] still works. We carry out a simulation study to compare the performance of the various procedures.

We also show that the working model of a piecewise constant function with a single jump can be naturally extended to allow a smooth parametric curve to the left of the jump and a smooth parametric curve to the right of the jump. A model of this type is a two-phase linear regression (also called break-point regression), as has been found useful, for example, in change-point analysis for climate data [12] and the estimation of mixed layer depth from oceanic profile data [20]. Similar models are used in econometrics, where they are called structural change models and threshold regression models.

In change-point analysis the aim is to estimate the locations of jump discontinuities in an otherwise smooth curve. Methods to do this are well developed in the nonparametric regression literature; see, for example, [1, 5, 9]. No distinction is made between the working model that has the jump point and the model that is assumed to generate the data. In contrast, confidence intervals for split points are model-robust in the sense that they apply under misspecification of the discontinuous working model by a smooth curve. Split-point analysis can thus be seen as complementary to change-point analysis: it is more appropriate in applications (such as the Everglades example mentioned above) in which the regression function is thought to be smooth, and does not require the a priori existence of a jump discontinuity. The working model has the jump discontinuity and is simply designed to condense key information about the underlying curve to a small number of parameters.

Confidence intervals for change-points are highly unstable under model misspecification by a smooth curve due to a sharp decrease in estimator rate of convergence: from close to $n$ under the assumed change-point model, to only a cube-root rate under a smooth curve (as for split-point estimators). This is not surprising because the split point depends on local features of a smooth regression curve which are harder to estimate than jumps. Misspecification of a change-point model thus causes confidence intervals to be misleadingly narrow, and rules out applications in which the existence of an abrupt change cannot be assumed a priori. In contrast, misspecification of a continuous (parametric) regression model (e.g., linear regression) causes no change in the $\sqrt{n}$-rate of convergence and the model-robust (Huber–White) sandwich estimate of variance is available. While the statistical literature on change-point analysis and model-robust estimation is comprehensive, split-point estimation falls in the gap between these two topics and is in need of further development.

The paper is organized as follows. In Section 2 we develop our main results and indicate how they can be applied in generalized nonparametric regression settings. In Section 3 we discuss an extension of our procedures to decision trees that incorporate general parametric working models. Simulation results and an application to Everglades data are presented in Section 4. Proofs are collected in Section 5.

**2. Split-point estimation in nonparametric regression.** We start this section by studying the problem of estimating the split point in a binary decision tree for nonparametric regression.

Let $X, Y$ denote the (one-dimensional) predictor and response variables, respectively, and assume that $Y$ has a finite second moment. The nonparametric regression function $f(x) = E(Y|X = x)$ is to be approximated using a decision tree with a single (terminal) node, that is, a piecewise constant function with a single jump. The predictor $X$ is assumed to have a density $p_X$, and its distribution function is denoted $F_X$. For convenience, we adopt the usual representation $Y = f(X) + \varepsilon$, with the error $\varepsilon = Y - E(Y|X)$ having zero conditional mean given $X$. The conditional variance of $\varepsilon$ given $X = x$ is denoted $\sigma^2(x)$.

Suppose we have $n$ i.i.d. observations $(X_1, Y_1), (X_2, Y_2), \ldots, (X_n, Y_n)$ of $(X, Y)$. Consider the working model in which $f$ is treated as a stump, that is, a piecewise constant function with a single jump, having parameters $(\beta_l, \beta_u, d)$, where $d$ is the point at which the function jumps, $\beta_l$ is the value to the left of the jump and $\beta_u$ is the value to the right of the jump. Best projected values are then defined by

$$(2.1) \qquad (\beta_l^0, \beta_u^0, d^0) = \arg\min_{\beta_l, \beta_u, d} E[Y - \beta_l 1(X \leq d) - \beta_u 1(X > d)]^2.$$

Before proceeding, we impose some mild conditions.

CONDITIONS.

(A1) There is a unique minimizer $(\beta_l^0, \beta_u^0, d^0)$ of the expectation on the right-hand side of (2.1) with $\beta_l^0 \neq \beta_u^0$.

(A2) $f(x)$ is continuous and is continuously differentiable in an open neighborhood $N$ of $d^0$. Also, $f'(d^0) \neq 0$.

(A3) $p_X(x)$ does not vanish and is continuously differentiable on $N$.

(A4) $\sigma^2(x)$ is continuous on $N$.

(A5) $\sup_{x \in N} E[\varepsilon^2 1\{|\varepsilon| > \eta\}|X = x] \to 0$ as $\eta \to \infty$.

The vector $(\beta_l^0, \beta_u^0, d^0)$ then satisfies the normal equations

$$\beta_l^0 = E(Y|X \leq d^0), \qquad \beta_u^0 = E(Y|X > d^0), \qquad f(d^0) = \frac{\beta_l^0 + \beta_u^0}{2}.$$

The usual estimates of these quantities are obtained via least squares as

$$(2.2) \qquad (\hat{\beta}_l, \hat{\beta}_u, \hat{d}_n) = \arg\min_{\beta_l, \beta_u, d} \sum_{i=1}^n [Y_i - \beta_l 1(X_i \leq d) - \beta_u 1(X_i > d)]^2.$$

Here and in the sequel, whenever we refer to a minimizer, we mean some choice of minimizer rather than the set of all minimizers (similarly for maximizers). Our first result gives the joint asymptotic distribution of these least squares estimators.

THEOREM 2.1. *If* (A1)–(A5) *hold, then*

$$n^{1/3}(\hat{\beta}_l - \beta_l^0, \hat{\beta}_u - \beta_u^0, \hat{d}_n - d^0) \xrightarrow{d} (c_1, c_2, 1) \arg\max_t Q(t),$$

*where*

$$Q(t) = aW(t) - bt^2,$$

$W$ is a standard two-sided Brownian motion process on the real line, $a^2 = \sigma^2(d^0)p_X(d^0)$,

$$b = b_0 - \frac{1}{8}|\beta_l^0 - \beta_u^0|p_X(d^0)^2\left(\frac{1}{F_X(d^0)} + \frac{1}{1-F_X(d^0)}\right) > 0,$$

with $b_0 = |f'(d^0)|p_X(d^0)/2$ and

$$c_1 = \frac{p_X(d^0)(\beta_u^0 - \beta_l^0)}{2F_X(d^0)}, \qquad c_2 = \frac{p_X(d^0)(\beta_u^0 - \beta_l^0)}{2(1-F_X(d^0))}.$$

In our notation, Bühlmann and Yu's [3] Theorem 3.1 states that $n^{1/3}(\hat{d}_n - d^0) \to_d \arg\max_t Q_0(t)$, where $Q_0(t) = aW(t) - b_0 t^2$. The first step in their proof assumes that it suffices to study the case in which $(\beta_l^0, \beta_u^0)$ is known. To justify this, they claim that $(\hat{\beta}_l, \hat{\beta}_u)$ converges at $\sqrt{n}$-rate to the population projected values $(\beta_l^0, \beta_u^0)$, which is faster than the $n^{1/3}$-rate of convergence of $\hat{d}_n$ to $d^0$. However, Theorem 2.1 shows that this is not the case; all three parameter estimates converge at cube-root rate, and have a nondegenerate joint asymptotic distribution concentrated on a line through the origin. Moreover, the limiting distribution of $\hat{d}_n$ differs from the one stated by Bühlmann and Yu because $b \neq b_0$; their limiting distribution will appear later in connection with (2.8).

*Wald-type confidence intervals.* It can be shown using Brownian scaling (see, e.g., [2]) that

(2.3) $$Q(t) \stackrel{d}{=} a(a/b)^{1/3}Q_1((b/a)^{2/3}t),$$

where $Q_1(t) = W(t) - t^2$, so the limit in the above theorem can be expressed more simply as

$$(c_1, c_2, 1)(a/b)^{2/3}\arg\max_t Q_1(t).$$

Let $p_{\alpha/2}$ denote the upper $\alpha/2$-quantile of the distribution of $\arg\max_t Q_1(t)$ (this is symmetric about 0), known as Chernoff's distribution. Accurate values of $p_{\alpha/2}$, for selected values of $\alpha$, are available in [10], where numerical aspects of Chernoff's distribution are studied. Utilizing the above theorem, this allows us to construct approximate $100(1-\alpha)\%$ confidence limits simultaneously for all the parameters $(\beta_l^0, \beta_u^0, d^0)$ in the working model:

(2.4)
$$\hat{\beta}_l \pm \hat{c}_1\hat{\delta}_n, \qquad \hat{\beta}_u \pm \hat{c}_2\hat{\delta}_n, \qquad \hat{d}_n \pm \hat{\delta}_n,$$
$$\text{where } \hat{\delta}_n = n^{-1/3}(\hat{a}/\hat{b})^{2/3}p_{\alpha/2},$$

given consistent estimators $\hat{c}_1, \hat{c}_2, \hat{a}, \hat{b}$ of the nuisance parameters. The density and distribution function of $X$ at $d^0$ can be estimated without difficulty, since an i.i.d. sample from the distribution of $X$ is available. The derivative $f'(d^0)$ and the conditional variance $\sigma^2(d^0)$ are harder to estimate, but many methods to do this are available in the literature, for example, local polynomial fitting with data-driven local bandwidth selection [18].

These confidence intervals are centered on the point estimate and have the disadvantage of not adapting to any skewness in the sampling distribution, which might be a problem in small samples. A more serious problem, however, is that the width of the interval is proportional to $\hat{a}/\hat{b}$, which blows up if $\hat{b}$ is small relative to $\hat{a}$. It follows from Theorem 2.1 that in the presence of conditions (A2)–(A5), the uniqueness condition (A1) fails if $b < 0$. Moreover, $b < 0$ if the gradient of the regression function is less than the jump in the working model multiplied by the density of $X$ at the split point, $|f'(d^0)| < p_X(d^0)|\beta_u^0 - \beta_l^0|$. This suggests that the Wald-type confidence interval becomes unstable if the regression function is flat enough at the split point.

*Subsampling.* Theorem 2.1 also makes it possible to avoid the estimation of nuisance parameters by using the subsampling bootstrap, which involves drawing a large number of subsamples of size $m = m_n$ from the original sample of size $n$ (without replacement). Then we can estimate the limiting quantiles of $n^{1/3}(\hat{d}_n - d^0)$ using the empirical distribution of $m^{1/3}(\hat{d}_m^* - \hat{d}_n)$; here $\hat{d}_m^*$ is the value of the split-point of the best fitting stump based on the subsample. For consistent estimation of the quantiles, we need $m/n \to 0$. In the literature $m$ is referred to as the block-size; see [15]. The choice of $m$ has a strong effect on the precision of the confidence interval, so a data-driven choice of $m$ is recommended in practice; Delgado, Rodríguez-Poo and Wolf [4] suggest a bootstrap-based algorithm for this purpose.

*Confidence sets based on residual sums of squares.* Another strategy is to use the quadratic loss function as an asymptotic pivot, which can be inverted to provide a confidence set. Such an approach was originally suggested by Stein [19] for a multivariate normal mean and has recently been used by Genovese and Wasserman [8] for nonparametric wavelet regression. To motivate the approach in the present setting, consider testing the null hypothesis that the working model parameters take the values $(\beta_l, \beta_u, d)$. Under the working model with a constant error variance, the likelihood-ratio statistic for testing this null hypothesis is given by

$$\text{RSS}_0(\beta_l, \beta_u, d) = \sum_{i=1}^n (Y_i - \beta_l 1(X_i \leq d) - \beta_u 1(X_i > d))^2$$
$$- \sum_{i=1}^n (Y_i - \hat{\beta}_l 1(X_i \leq \hat{d}_n) - \hat{\beta}_u 1(X_i > \hat{d}_n))^2.$$

The corresponding profiled RSS statistic for testing the null hypothesis that $d^0 = d$ replaces $\beta_l$ and $\beta_u$ in $\text{RSS}_0$ by their least squares estimates under the null hypothesis, giving

$$\text{RSS}_1(d) = \sum_{i=1}^{n}(Y_i - \hat{\beta}_l^d 1(X_i \leq d) - \hat{\beta}_u^d 1(X_i > d))^2$$
$$- \sum_{i=1}^{n}(Y_i - \hat{\beta}_l 1(X_i \leq \hat{d}_n) - \hat{\beta}_u 1(X_i > \hat{d}_n))^2,$$

where

$$(\hat{\beta}_l^d, \hat{\beta}_u^d) = \arg\min_{\beta_l, \beta_u} \sum_{i=1}^{n}(Y_i - \beta_l 1(X_i \leq d) - \beta_u 1(X_i > d))^2.$$

Our next result provides the asymptotic distribution of these residual sums of squares.

THEOREM 2.2. *If* (A1)–(A5) *hold, then*

$$n^{-1/3}\text{RSS}_0(\beta_l^0, \beta_u^0, d^0) \xrightarrow{d} 2|\beta_l^0 - \beta_u^0|\max_t Q(t),$$

*where $Q$ is given in Theorem* 2.1, *and $n^{-1/3}\text{RSS}_1(d^0)$ has the same limiting distribution.*

Using the Brownian scaling (2.3), the above limiting distribution can be expressed more simply as

$$2|\beta_l^0 - \beta_u^0|a(a/b)^{1/3}\max_t Q_1(t).$$

This leads to the following approximate $100(1-\alpha)\%$ confidence set for the split-point:

(2.5) $$\{d : \text{RSS}_1(d) \leq 2n^{1/3}|\hat{\beta}_l - \hat{\beta}_u|\hat{a}(\hat{a}/\hat{b})^{1/3}q_\alpha\},$$

where $q_\alpha$ is the upper $\alpha$-quantile of $\max_t Q_1(t)$. This confidence set becomes unstable if $\hat{b}$ is small relative to $\hat{a}$, as with the Wald-type confidence interval. This problem can be lessened by changing the second term in $\text{RSS}_1$ to make use of the information in the null hypothesis, to obtain

$$\text{RSS}_2(d) = \sum_{i=1}^{n}(Y_i - \hat{\beta}_l^d 1(X_i \leq d) - \hat{\beta}_u^d 1(X_i > d))^2$$
$$- \sum_{i=1}^{n}(Y_i - \hat{\beta}_l^d 1(X_i \leq \hat{d}_n^d) - \hat{\beta}_u^d 1(X_i > \hat{d}_n^d))^2,$$

where

$$\hat{d}_n^d = \arg\min_{d'} \sum_{i=1}^n (Y_i - \hat{\beta}_l^d 1(X_i \leq d') - \hat{\beta}_u^d 1(X_i > d'))^2. \tag{2.6}$$

The following result gives the asymptotic distribution of $\text{RSS}_2(d^0)$.

THEOREM 2.3. *If* (A1)–(A5) *hold, then*

$$n^{-1/3}\text{RSS}_2(d^0) \xrightarrow{d} 2|\beta_l^0 - \beta_u^0| \max_t Q_0(t),$$

*where* $Q_0(t) = aW(t) - b_0 t^2$, *and* $a$, $b_0$ *are given in Theorem* 2.1.

This leads to the following approximate $100(1-\alpha)\%$ confidence set for the split point:

$$\{d : \text{RSS}_2(d) \leq 2n^{1/3}|\hat{\beta}_l - \hat{\beta}_u|\hat{a}(\hat{a}/\hat{b}_0)^{1/3}q_\alpha\}, \tag{2.7}$$

where $\hat{b}_0$ is a consistent estimator of $b_0$. This confidence set could be unstable if $\hat{b}_0$ is small compared with $\hat{a}$, but this is less likely to occur than the instability we described earlier because $b_0 > b$. The proof of Theorem 2.3 also shows that $n^{1/3}(\hat{d}_n^{d^0} - d^0)$ converges in distribution to $\arg\max_t Q_0(t)$, recovering the limit distribution in Theorem 3.1 of [3], and this provides another pivot-type confidence set for the split point,

$$\{d : |\hat{d}_n^d - d| \leq n^{-1/3}(\hat{a}/\hat{b}_0)^{2/3} p_{\alpha/2}\}. \tag{2.8}$$

Typically, (2.5), (2.7) and (2.8) are not intervals, but their endpoints, or the endpoints of their largest component, can be used as approximate confidence limits.

REMARK 1. The uniqueness condition (A1) may be violated if the regression function is not monotonic on the support of $X$. A simple example in which uniqueness fails is given by $f(x) = x^2$ and $X \sim \text{Unif}[-1,1]$, in which case the normal equations for the split point have two solutions, $d^0 = \pm 1/\sqrt{2}$, and the corresponding $\beta_l^0$ and $\beta_u^0$ are different for each solution; neither split point has a natural interpretation because the regression function has no trend. More generally, we would expect lack of unique split points for regression functions that are unimodal on the interior of the support of $X$. In a practical situation, split-point analysis (with stumps) should not be used unless there is reason to believe that a trend is present, in which case we expect there to be a unique split point. An increasing trend, for instance, gives that $E(Y|X \leq d) < E(Y|X > d)$ for all $d$, so a unique split point will exist provided the normal equation $g(d) = 0$ has a unique solution, where $g$ is the "centered" regression function $g(d) = f(d) - (E(Y|X \leq d) + E(Y|X > d))/2$. A sufficient condition for $g(d) = 0$ to have a unique solution is that $g$ is continuous and strictly increasing, with $g(x_0) < 0$ and $g(x_1) > 0$ for some $x_0 < x_1$ in the support of $X$.

*Generalized nonparametric regression.* Our results apply to split point estimation for a generalized nonparametric regression model in which the conditional distribution of $Y$ given $X$ is assumed to belong to an exponential family. The canonical parameter of the exponential family is expressed as $\theta(X)$ for an unknown smooth function $\theta(\cdot)$, and we are interested in estimation of the split point in a decision tree approximation of $\theta(\cdot)$. Nonparametric estimation of $\theta(\cdot)$ has been studied extensively; see, for example, [6], Section 5.4. Important examples include the binary choice or nonlinear logistic regression model $Y|X \sim \text{Ber}(f(X))$, where $f(x) = e^{\theta(x)}/(1+e^{\theta(x)})$, and the Poisson regression model $Y|X \sim \text{Poi}(f(X))$, where $f(x) = e^{\theta(x)}$.

The conditional density of $Y$ given $X = x$ is specified as

$$p(y|x) = \exp\{\theta(x)y - B(\theta(x))\}h(y),$$

where $B(\cdot)$ and $h(\cdot)$ are known functions. Here $p(\cdot|x)$ is a probability density function with respect to some given Borel measure $\mu$. Here the cumulant function $B$ is twice continuously differentiable and $B'$ is strictly increasing on the range of $\theta(\cdot)$. It can be shown that $f(x) = E(Y|X=x) = B'(\theta(x))$, or equivalently $\theta(x) = \psi(f(x))$, where $\psi = (B')^{-1}$ is the link function. For logistic regression $\psi(t) = \log(t/(1-t))$ is the logit function, and for Poisson regression $\psi(t) = \log(t)$. The link function is known, continuous and strictly increasing, so a stump approximation to $\theta(x)$ is equivalent to a stump approximation to $f(x)$, and the split points are identical. Exploiting this equivalence, we define the best projected values of the stump approximation for $\theta(\cdot)$ as $(\psi(\beta_l^0), \psi(\beta_u^0), d^0)$, where $(\beta_l^0, \beta_u^0, d^0)$ are given in (2.1).

Our earlier results apply under a reduced set of conditions due to the additional structure in the exponential family model: we only need (A1), (A2) with $\theta(\cdot)$ in place of $f$, and (A3). It is then easy to check that the original assumption (A2) holds; in particular, $f'(d^0) = B''(\theta(d^0))\theta'(d^0) \neq 0$. To check (A4), note that $\sigma^2(x) = \text{Var}(Y|X=x) = B''(\theta(x))$ is continuous in $x$. Finally, to check (A5), let $N$ be a bounded neighborhood of $d^0$. Note that $f(\cdot)$ and $\theta(\cdot)$ are bounded on $N$. Let $\theta_0 = \inf_{x \in N} \theta(x)$ and $\theta_1 = \sup_{x \in N} \theta(x)$. For $\eta$ sufficiently large, $\{y : |y - f(x)| > \eta\} \subset \{y : |y| > \eta/2\}$ for all $x \in N$, and consequently

$$\sup_{x \in N} E[\varepsilon^2 1\{|\varepsilon| > \eta\}|X=x] = \sup_{x \in N} \int_{|y-f(x)|>\eta} (y - f(x))^2 p(y|x) \, d\mu(y)$$
$$\leq C \int_{|y|>\eta/2} (y^2 + 1)(e^{\theta_0 y} + e^{\theta_1 y})h(y) \, d\mu(y) \to 0$$

as $\eta \to \infty$, where $C$ is a constant (not depending on $\eta$). The last step follows from the dominated convergence theorem.

We have focused on confidence sets for the split point, but $\beta_l^0$ and $\beta_u^0$ may also be important. For example, in logistic regression where the response $Y$ is an indicator variable, the relative risk

$$r = P(Y=1|X>d^0)/P(Y=1|X\le d^0) = \beta_u^0/\beta_l^0$$

is useful for comparing the risks before and after the split point. Using Theorem 2.1 and the delta method, we can obtain the approximate $100(1-\alpha)\%$ confidence limits

$$\exp\left(\log(\hat\beta_u/\hat\beta_l) \pm \left(\frac{\hat c_2}{\hat\beta_u} - \frac{\hat c_1}{\hat\beta_l}\right)\hat\delta_n\right)$$

for $r$, where $\hat\delta_n$ is defined in (2.4) and it is assumed that $c_1/\beta_l \ne c_2/\beta_u$ to ensure that $\hat\beta_u/\hat\beta_l$ has a nondegenerate limit distribution. The odds ratio for comparing $P(Y=1|X\le d^0)$ and $P(Y=1|X>d^0)$ can be treated in a similar fashion.

**3. Extending the decision tree approach.** We have noted that split-point estimation with stumps should only be used if a trend is present. The split-point approach can be adapted to more complex situations, however, by using a more flexible working model that provides a better approximation to the underlying regression curve. In this section, we indicate how our main results extend to a broad class of parametric working models. The proofs are omitted as they run along similar lines.

The constants $\beta_l$ and $\beta_u$ are now replaced by functions $\Psi_l(\beta_l, x)$ and $\Psi_u(\beta_u, x)$ specified in terms of vector parameters $\beta_l$ and $\beta_u$. These functions are taken to be twice continuously differentiable with respect to $\beta_l \in \mathbb{R}^m$ and $\beta_u \in \mathbb{R}^k$, respectively, and continuously differentiable with respect to $x$. The best projected values of the parameters in the working model are defined by

$$\begin{aligned}(\beta_l^0, \beta_u^0, d^0) = \arg\min_{\beta_l,\beta_u,d} E[Y &- \Psi_l(\beta_l, X)\mathbf{1}(X\le d) \\ &- \Psi_u(\beta_u, X)\mathbf{1}(X>d)]^2,\end{aligned} \quad (3.1)$$

and the corresponding normal equations are

$$E\left[\frac{\partial}{\partial\beta_l}\Psi_l(\beta_l^0, X)(Y - \Psi_l(\beta_l^0, X))\mathbf{1}(X\le d^0)\right] = 0,$$

$$E\left[\frac{\partial}{\partial\beta_u}\Psi_u(\beta_u^0, X)(Y - \Psi_u(\beta_u^0, X))\mathbf{1}(X>d^0)\right] = 0,$$

and $f(d^0) = \Psi(d^0)$, where $\Psi(x) = (\Psi_l(\beta_l^0, x) + \Psi_u(\beta_u^0, x))/2$. The least squares estimates of these quantities are obtained as

$$\begin{aligned}(\hat\beta_l, \hat\beta_u, \hat d_n) = \arg\min_{\beta_l,\beta_u,d} \sum_{i=1}^n [Y_i &- \Psi_l(\beta_l, X_i)\mathbf{1}(X_i\le d) \\ &- \Psi_u(\beta_u, X_i)\mathbf{1}(X_i>d)]^2.\end{aligned} \quad (3.2)$$

To extend Theorem 2.1, we need to modify conditions (A1) and (A2) as follows:

(A1)′ There is a unique minimizer $(\beta_l^0, \beta_u^0, d^0)$ of the expectation on the right-hand side of (3.1) with $\Psi_l(\beta_l^0, d^0) \neq \Psi_u(\beta_u^0, d^0)$.

(A2)′ $f(x)$ is continuously differentiable in an open neighborhood $N$ of $d^0$. Also, $f'(d^0) \neq \Psi'(d^0)$.

In addition, we need the following Lipschitz condition on the working model:

(A6) There exist functions $\dot{\Psi}_l(x)$ and $\dot{\Psi}_u(x)$, bounded on compacts, such that

$$|\Psi_l(\beta_l, x) - \Psi_l(\tilde{\beta}_l, x)| \leq \dot{\Psi}_l(x)|\beta_l - \tilde{\beta}_l|$$

and

$$|\Psi_u(\beta_u, x) - \Psi_u(\tilde{\beta}_u, x)| \leq \dot{\Psi}_u(x)|\beta_u - \tilde{\beta}_u|$$

with $\dot{\Psi}_l(X), \Psi_l(\beta_l^0, X), \dot{\Psi}_u(X), \Psi_u(\beta_u^0, X)$ having finite fourth moments, where $|\cdot|$ is Euclidean distance.

Condition (A6) holds, for example, if $\Psi_l(\beta_l, x)$ and $\Psi_u(\beta_u, x)$ are polynomials in $x$ with the components of $\beta_l$ and $\beta_u$ serving as coefficients, and $X$ has a finite moment of sufficiently high order.

THEOREM 3.1. *If (A1)′, (A2)′ and (A3)–(A6) hold, then*

$$n^{1/3}(\hat{\beta}_l - \beta_l^0, \hat{\beta}_u - \beta_u^0, \hat{d}_n - d^0) \xrightarrow{d} \arg\min_h \tilde{W}(h),$$

*where $\tilde{W}$ is the Gaussian process*

$$\tilde{W}(h) = \tilde{a}W(h_{m+k+1}) + h^T V h/2, \qquad h \in \mathbb{R}^{m+k+1},$$

*$V$ is the (positive definite) Hessian matrix of the function*

$$(\beta_l, \beta_u, d) \mapsto E[Y - \Psi_l(\beta_l, X)1(X \leq d) - \Psi_u(\beta_u, X)1(X > d)]^2$$

*evaluated at $(\beta_l^0, \beta_u^0, d^0)$, and $\tilde{a} = 2|\Psi_l(\beta_l^0, d^0) - \Psi_u(\beta_u^0, d^0)|(\sigma^2(d^0)p_X(d^0))^{1/2}$.*

REMARK 2. As in the decision tree case, subsampling can now be used to construct confidence intervals for the parameters of the working model. Although Brownian scaling is still available [minimizing $\tilde{W}(h)$ by first holding $h_{m+k+1}$ fixed], the construction of Wald-type confidence intervals would be cumbersome, needing estimation of all the nuisance parameters involved in $\tilde{a}$ and $V$. The complexity of $V$ is already evident when $\beta_l$ and $\beta_u$ are one-dimensional, in which case direct computation shows that $V$ is the $3 \times 3$ matrix with entries $V_{12} = V_{21} = 0$,

$$V_{11} = 2\int_{-\infty}^{d^0} \left(\frac{\partial}{\partial \beta_l}\Psi_l(\beta_l^0, x)\right)^2 p_X(x)\,dx$$
$$+ 2\int_{-\infty}^{d^0} \frac{\partial^2}{\partial \beta_l^2}\Psi_l(\beta_l^0, x)(\Psi_l(\beta_l^0, x) - f(x))p_X(x)\,dx,$$
$$V_{22} = 2\int_{d^0}^{\infty} \left(\frac{\partial}{\partial \beta_u}\Psi_u(\beta_u^0, x)\right)^2 p_X(x)\,dx$$
$$+ 2\int_{d^0}^{\infty} \frac{\partial^2}{\partial \beta_u^2}\Psi_u(\beta_u^0, x)(\Psi_u(\beta_u^0, x) - f(x))p_X(x)\,dx,$$
$$V_{33} = 2|(\Psi_l(\beta_u^0, d^0) - \Psi_l(\beta_l^0, d^0))(f'(d^0) - \Psi'(d^0))|p_X(d^0),$$
$$V_{13} = V_{31} = (\Psi_l(\beta_l^0, d^0) - \Psi_u(\beta_u^0, d^0))\frac{\partial}{\partial \beta_l}\Psi_l(\beta_l^0, d^0)p_X(d^0),$$
$$V_{23} = V_{32} = (\Psi_l(\beta_l^0, d^0) - \Psi_u(\beta_u^0, d^0))\frac{\partial}{\partial \beta_u}\Psi_u(\beta_u^0, d^0)p_X(d^0).$$

Next we show that extending Theorem 2.3 allows us to circumvent this problem. Two more conditions are needed:

(A7) $\int_D (\Psi_l(\beta_l^0, x) - \Psi_u(\beta_u^0, x))(f(x) - \Psi(x))p_X(x)\,dx \neq 0$, for $D = (-\infty, d^0]$ and $D = [d^0, \infty)$.

(A8) $\sqrt{n}(\hat{\beta}_l^{d^0} - \beta_l^0) = O_p(1)$ and $\sqrt{n}(\hat{\beta}_u^{d^0} - \beta_u^0) = O_p(1)$, where $\hat{\beta}_l^d$ and $\hat{\beta}_u^d$ are defined in an analogous fashion to Section 2.

Note that (A8) holds automatically in the setting of Section 2, using the central limit theorem and the delta method. In the present setting, sufficient conditions for (A8) can be easily formulated in terms of $\Psi_l$, $\Psi_u$ and the joint distribution of $(X, Y)$, using the theory of $Z$-estimators. If we define $\phi_{\beta_l}(x, y) = (y - \Psi_l(\beta_l, x))(\partial \Psi_l(\beta_l, x)/\partial \beta_l)1(x \leq d^0)$, then $\beta_l^0$ satisfies the normal equation $P\phi_{\beta_l} = 0$, while $\hat{\beta}_l^{d^0}$ satisfies $\mathbb{P}_n \phi_{\beta_l} = 0$, where $\mathbb{P}_n$ is the empirical distribution of $(X_i, Y_i)$. Sufficient conditions for the asymptotic normality of $\sqrt{n}(\hat{\beta}_l^{d^0} - \beta_l^0)$ are then given by Lemma 3.3.5 of [21] (see also Examples 3.3.7 and 3.3.8 in Section 3.3 in [21], which are special cases of Lemma 3.3.5 in the context of finite-dimensional parametric models) in conjunction with $\beta \mapsto P\phi_\beta$ possessing a nonsingular derivative at $\beta_l^0$. In particular, if $\Psi_l$ and $\Psi_u$ are polynomials in $x$ with the $\beta_l$ and $\beta_u$ serving as coefficients, then the displayed condition in Example 3.3.7 is easily verifiable under the assumption that $X$ has a finite moment of a sufficiently high order (which is trivially true if $X$ has compact support).

Defining

$$\text{RSS}_2(d) = \sum_{i=1}^{n}(Y_i - \Psi_l(\hat{\beta}_l^d, X_i)\mathbf{1}(X_i \leq d) - \Psi_u(\hat{\beta}_u^d, X_i)\mathbf{1}(X_i > d))^2$$

$$- \sum_{i=1}^{n}(Y_i - \Psi_l(\hat{\beta}_l^d, X_i)\mathbf{1}(X_i \leq \hat{d}_n^d) - \Psi_u(\hat{\beta}_u^d, X_i)\mathbf{1}(X_i > \hat{d}_n^d))^2,$$

where

$$\hat{d}_n^d = \arg\min_{d'} \sum_{i=1}^{n}(Y_i - \Psi_l(\hat{\beta}_l^d, X_i)\mathbf{1}(X_i \leq d') - \Psi_u(\hat{\beta}_u^d, X_i)\mathbf{1}(X_i > d'))^2,$$

we obtain the following extension of Theorem 2.3.

THEOREM 3.2. *If (A1)′, (A2)′, (A3)–(A5), (A7) and (A8) hold, and the random variables $\dot{\Psi}_l(X)$, $\Psi_l(\beta_l^0, X)$, $\dot{\Psi}_u(X)$ and $\Psi_u(\beta_u^0, X)$ are square integrable, then*

$$n^{-1/3}\text{RSS}_2(d^0) \xrightarrow{d} 2|\Psi_l(\beta_l^0, d^0) - \Psi_u(\beta_u^0, d^0)|\max_t Q_0(t),$$

*where $Q_0(t) = aW(t) - b_0 t^2$, and $a^2 = \sigma^2(d^0)p_X(d^0)$, $b_0 = |f'(d^0) - \Psi'(d^0)| \times p_X(d^0)$.*

Application of the above result to construct confidence sets [as in (2.7)] is easier than using Theorem 3.1, since estimation of $a$ and $b_0$ requires much less work than estimation of the matrix $V$; the latter is essentially intractable, even for moderate $k$ and $m$.

**4. Numerical examples.** In this section we compare the various confidence sets for the split point in a binary decision tree using simulated data. We also develop the Everglades application mentioned in the Introduction.

4.1. *Simulation study.* We consider a regression model of the form $Y = f(X) + \varepsilon$, where $X \sim \text{Unif}[0,1]$ and $\varepsilon | X \sim N(0, \sigma^2(X))$. The regression function $f$ is specified as the sigmoid (or logistic distribution) function

$$f(x) = e^{15(x-0.5)}/(1 + e^{15(x-0.5)}).$$

This increasing S-shaped function rises steeply between 0.2 and 0.8, but is relatively flat otherwise. It is easily checked that $d^0 = 0.5$, $\beta_l^0 = 0.092$ and $\beta_u^0 = 0.908$. We take $\sigma^2(x) = 0.25$ to produce an example with homoscedastic error, and $\sigma^2(x) = \exp(-2.77x)$ for an example with heteroscedastic error; these two error variances agree at the split point.

TABLE 1
Coverage and average confidence interval length, $\sigma^2(x) = 0.25$

| | Subsampling | | Wald | | $RSS_1$ | | $RSS_2$ | |
|---|---|---|---|---|---|---|---|---|
| $n$ | Coverage | Length | Coverage | Length | Coverage | Length | Coverage | Length |
| 75 | 0.957 | 0.326 | 0.883 | 0.231 | 0.942 | 0.273 | 0.957 | 0.345 |
| 100 | 0.970 | 0.283 | 0.894 | 0.210 | 0.954 | 0.235 | 0.956 | 0.280 |
| 200 | 0.978 | 0.200 | 0.926 | 0.167 | 0.952 | 0.174 | 0.959 | 0.198 |
| 500 | 0.991 | 0.136 | 0.947 | 0.123 | 0.947 | 0.118 | 0.948 | 0.128 |
| 1000 | 0.929 | 0.093 | 0.944 | 0.097 | 0.955 | 0.091 | 0.952 | 0.098 |
| 1500 | 0.936 | 0.098 | 0.947 | 0.085 | 0.933 | 0.078 | 0.921 | 0.083 |
| 2000 | 0.944 | 0.090 | 0.954 | 0.077 | 0.935 | 0.070 | 0.939 | 0.074 |

To compute the subsampling confidence interval, a data-driven choice of block size was not feasible computationally. Instead, the block size was determined via a pilot simulation. For a given sample size, 1000 independently replicated samples were generated from the (true) regression model, and for each data set a collection of subsampling-based intervals (of nominal level 95%) was constructed, for block sizes of the form $m_n = n^\gamma$, for $\gamma$ on a grid of values between 0.33 and 0.9. The block size giving the greatest empirical accuracy (in terms of being closest to 95% coverage based on the replicated samples) was used in the subsequent simulation study. To provide a fair comparison, we used the true values of the nuisance parameters to calibrate the Wald- and RSS-type confidence sets. For $RSS_1$ and $RSS_2$ we use the endpoints of the longest connected component to specify confidence limits.

Tables 1 and 2 report the results of simulations based on 1000 replicated samples, with sample sizes ranging from 75 to 2000, and each confidence interval (CI) calibrated to have nominal 95% coverage. The subsampling CI tends to be wider than the others, especially at small sample sizes. The Wald-type CI suffers from severe undercoverage, especially in the heteroscedastic case and at small sample sizes. The $RSS_1$-type CI is also prone to undercoverage in the heteroscedastic case. The $RSS_2$-type CI performs well, although there is a slight undercoverage at high sample sizes (the interval formed by the endpoints of the entire confidence set has greater accuracy in that case).

4.2. *Application to Everglades data.* The "river of grass" known as the Everglades is a majestic wetland covering much of South Florida. Severe damage to large swaths of this unique ecosystem has been caused by pollution from agricultural fertilizers and the disruption of water flow (e.g., from the construction of canals). Efforts to restore the Everglades started in earnest in the early 1990s. In 1994, the Florida legislature passed the Everglades Forever Act, which called for a threshold level of total phosphorus that would prevent an "imbalance in natural populations of aquatic flora or fauna." This threshold may eventually be set at around 10 or 15 parts per billion (ppb), but it remains undecided despite extensive scientific study and much political and legal debate; see [17] for a discussion of the statistical issues involved.

Between 1992 and 1998, the Duke University Wetlands Center (DUWC) carried out a dosing experiment at two unimpacted sites in the Everglades. This experiment was designed to find the threshold level of total phosphorus concentration at which biological imbalance occurs. Changes in the abundance of various phosphorus-sensitive species were monitored along dosing channels in which a gradient of phosphorus concentration had been established. Qian, King and Richardson [16] analyzed data from this experiment using Bayesian change-point analysis, and also split point estimation with the split point being interpreted as the threshold level at which biological imbalance occurs. Uncertainty in the split point was evaluated using Efron's bootstrap.

We illustrate our approach with one particular species monitored in the DUWC dosing experiment: the bladderwort *Utricularia Purpurea*, which is considered a keystone species for the health of the Everglades ecosystem. Figure 1 shows 340 observations of stem density plotted against the six-month geometric mean of total phosphorus concentration. The displayed data were collected in August 1995, March 1996, April 1998 and August 1998 (observations taken at unusually low or high water levels, or before the system stabilized in 1995, are excluded). Water levels fluctuate greatly and have a strong influence on species abundance, so a separate analysis for each data collection period would be preferable, but not enough data are available for separate analyses and a more sophisticated model would be needed, so for simplicity we have pooled all the data.

Estimates of $p_X$, $f'$ and $\sigma^2$ needed for $\hat{a}$, $\hat{b}$ and $\hat{b}_0$, and the estimate of $f$ shown in Figure 1, are found using David Ruppert's (Matlab) implemen-tation of local polynomial regression and density estimation with empirical-bias bandwidth selection [18]. The estimated regression function shows a fairly steady decrease in stem density with increasing phosphorus concentration, but there is no abrupt change around the split point estimate of 12.8 ppb, so we expect the CIs to be relatively wide. The 95% Wald-type and $RSS_1$-type CIs for the split point are 0.7–24.9 and 9.7–37.1 ppb, respectively. The instability problem mentioned earlier may be causing these CIs to be so wide (here $\hat{a}/\hat{b} = 722$). The subsampling and $RSS_2$-type CIs are narrower, at 8.5–17.1 and 7.1–26.1 ppb, respectively (see the vertical lines in Figure 1), but they still leave considerable uncertainty about the true location of the split point. The 10-ppb threshold recommended by the Florida Department of Environmental Protection [13] falls into these CIs.

TABLE 2
*Coverage and average confidence interval length, $\sigma^2(x) = \exp(-2.77x)$*

|  | Subsampling | | Wald | | $RSS_1$ | | $RSS_2$ | |
|---|---|---|---|---|---|---|---|---|
| $n$ | Coverage | Length | Coverage | Length | Coverage | Length | Coverage | Length |
| 75 | 0.951 | 0.488 | 0.863 | 0.231 | 0.929 | 0.270 | 0.949 | 0.354 |
| 100 | 0.957 | 0.315 | 0.884 | 0.210 | 0.923 | 0.231 | 0.944 | 0.283 |
| 200 | 0.977 | 0.257 | 0.915 | 0.167 | 0.939 | 0.173 | 0.949 | 0.196 |
| 500 | 0.931 | 0.124 | 0.926 | 0.123 | 0.936 | 0.117 | 0.948 | 0.128 |
| 1000 | 0.917 | 0.095 | 0.941 | 0.097 | 0.948 | 0.090 | 0.945 | 0.097 |
| 1500 | 0.938 | 0.083 | 0.938 | 0.085 | 0.928 | 0.078 | 0.922 | 0.083 |
| 2000 | 0.945 | 0.076 | 0.930 | 0.077 | 0.933 | 0.070 | 0.934 | 0.074 |

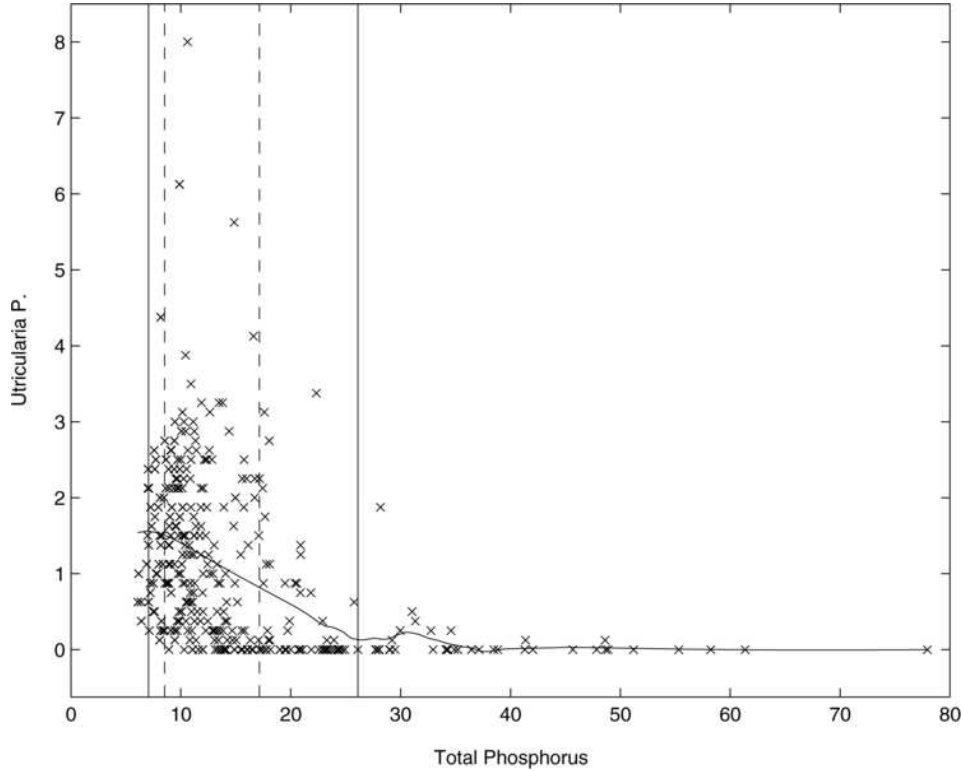

Fig. 1. *Data from the DUWC Everglades phosphorus dosing study showing variations in bladderwort (Utricularia P.) stem density (number of stems per square meter) in response to total phosphorus concentration (six-month geometric mean, units of ppb). The vertical solid lines show the limits of the $RSS_2$-type 95% confidence interval for the split point. The vertical dashed lines show the limits of the subsampling confidence interval. The local polynomial regression fit is also plotted.*

CONFIDENCE SETS FOR SPLIT POINTS 17The interpretation of the split point as a biological threshold is the source of some controversy in the debate over a numeric phosphorus criterion [13]. It can be argued that the split point is only crudely related to biological response and that it is a statistical construct depending on an artificial working model. Yet the split point approach fulfills a clear need in the absence of better biological understanding, and is preferable to a change-point analysis in this application, as discussed in the Introduction.

**5. Proofs.** The proofs have certain points in common with Bühlmann and Yu [3] and Kim and Pollard [11], but to make them more self-contained we mainly appeal to general results on empirical processes and $M$-estimation that are collected in [21].

We begin by proving Theorem 2.3, which is closely related to Theorem 3.1 of [3].

PROOF OF THEOREM 2.3. We derive the joint limiting distribution of

$$(n^{1/3}(\hat{d}_n^{d^0} - d^0), n^{-1/3}\mathrm{RSS}_2(d^0)),$$

the marginals of which are involved in calibrating the confidence sets (2.7) and (2.8). To simplify the notation, we denote $(\hat{\beta}_l^{d^0}, \hat{\beta}_u^{d^0}, \hat{d}_n^{d^0})$ by $(\hat{\beta}_l^0, \hat{\beta}_u^0, \hat{d}_n^0)$. Also, we assume that $\beta_l^0 > \beta_u^0$; the derivation for the other case is analogous. Letting $\mathbb{P}_n$ denote the empirical measure of the pairs $\{(X_i, Y_i), i = 1, \ldots, n\}$, we can write

$$\begin{aligned}
\mathrm{RSS}_2(d^0) &= \sum_{i=1}^n (Y_i - \hat{\beta}_l^0)^2 (1(X_i \leq d^0) - 1(X_i \leq \hat{d}_n^0)) \\
&\quad + \sum_{i=1}^n (Y_i - \hat{\beta}_u^0)^2 (1(X_i > d^0) - 1(X_i > \hat{d}_n^0)) \\
&= n\mathbb{P}_n[((Y - \hat{\beta}_l^0)^2 - (Y - \hat{\beta}_u^0)^2)(1(X \leq d^0) - 1(X \leq \hat{d}_n^0))] \\
&= 2(\hat{\beta}_l^0 - \hat{\beta}_u^0) n\mathbb{P}_n\left[\left(Y - \frac{\hat{\beta}_l^0 + \hat{\beta}_u^0}{2}\right)(1(X \leq \hat{d}_n^0) - 1(X \leq d^0))\right].
\end{aligned}$$

Therefore,

$$\begin{aligned}
&n^{-1/3}\mathrm{RSS}_2(d^0) \\
&\quad = 2(\hat{\beta}_l^0 - \hat{\beta}_u^0) n^{2/3} \mathbb{P}_n[(Y - \hat{f}(d^0))(1(X \leq \hat{d}_n^0) - 1(X \leq d^0))],
\end{aligned}$$

where $\hat{f}(d^0) = (\hat{\beta}_l^0 + \hat{\beta}_u^0)/2$. Let

$$\xi_n(d) = n^{2/3} \mathbb{P}_n[(Y - \hat{f}(d^0))(1(X \leq d) - 1(X \leq d^0))]$$

and let $\tilde{d}_n$ be the maximizer of this process. Since $\hat{\beta}_l^0 - \hat{\beta}_u^0 \to \beta_l^0 - \beta_u^0 > 0$ almost surely, it is easy to see that $\tilde{d}_n = \hat{d}_n^0$ for $n$ sufficiently large almost surely. Hence, the limiting distribution of $n^{-1/3}\text{RSS}_2(d^0)$ must be the same as that of $2(\hat{\beta}_l^0 - \hat{\beta}_u^0)\xi_n(\tilde{d}_n)$, which in turn is the same as that of $2(\beta_l^0 - \beta_u^0)\xi_n(\tilde{d}_n)$ [provided $\xi_n(\tilde{d}_n)$ has a limit distribution], because $\hat{\beta}_l^0$ and $\hat{\beta}_u^0$ are $\sqrt{n}$-consistent. Furthermore, the limiting distribution of $n^{1/3}(\hat{d}_n^0 - d^0)$ is the same as that of $n^{1/3}(\tilde{d}_n - d^0)$ (provided a limiting distribution exists).

Let $Q_n(t) = \xi_n(d^0 + tn^{-1/3})$ and $\hat{t}_n = \arg\max_t Q_n(t)$, so that $\hat{t}_n = n^{1/3}(\tilde{d}_n - d^0)$. It now suffices to find the joint limiting distribution of $(\hat{t}_n, Q_n(\hat{t}_n))$. Lemma 5.1 below shows that $Q_n(t)$ converges in distribution in the space $B_{\text{loc}}(\mathbb{R})$ (the space of locally bounded functions on $\mathbb{R}$ equipped with the topology of uniform convergence on compacta) to the Gaussian process $Q_0(t) \equiv aW(t) - b_0 t^2$ whose distribution is a tight Borel measure concentrated on $C_{\max}(\mathbb{R})$ [the separable subspace of $B_{\text{loc}}(\mathbb{R})$ of all continuous functions on $\mathbb{R}$ that diverge to $-\infty$ as the argument runs off to $\pm\infty$ and that have a unique maximum]. Furthermore, Lemma 5.1 shows that the sequence $\{\hat{t}_n\}$ of maximizers of $\{Q_n(t)\}$ is $O_p(1)$. By Theorem 5.1 below, we conclude that $(\hat{t}_n, Q_n(\hat{t}_n)) \to_d (\arg\max_t Q_0(t), \max_t Q_0(t))$. This completes the proof. □

The following theorem provides sufficient conditions for the joint weak convergence of a sequence of maximizers and the corresponding maxima of a general sequence of processes in $B_{\text{loc}}(\mathbb{R})$. A referee suggested that an alternative approach would be to use $D(\mathbb{R})$ (the space of right-continuous functions with left-limits equipped with Lindvall's extension of the Skorohod topology) instead of $B_{\text{loc}}(\mathbb{R})$, as in an argmax-continuous mapping theorem due to Ferger ([7], Theorem 3).

THEOREM 5.1. *Let $\{Q_n(t)\}$ be a sequence of stochastic processes converging in distribution in the space $B_{\text{loc}}(\mathbb{R}^k)$ to the process $Q(t)$, whose distribution is a tight Borel measure concentrated on $C_{\max}(\mathbb{R}^k)$. If $\{\hat{t}_n\}$ is a sequence of maximizers of $\{Q_n(t)\}$ such that $\hat{t}_n = O_p(1)$, then*

$$(\hat{t}_n, Q_n(\hat{t}_n)) \xrightarrow{d} \left(\arg\max_t Q(t), \max_t Q(t)\right).$$

PROOF. For simplicity, we provide the proof for the case that $k = 1$; the same argument essentially carries over to the $k$-dimensional case. By invoking Dudley's representation theorem (Theorem 2.2 of [11]), for the processes $Q_n$, we can construct a sequence of processes $\tilde{Q}_n$ and a process $\tilde{Q}$ defined on a common probability space $(\tilde{\Omega}, \tilde{\mathcal{A}}, \tilde{P})$ with (a) $\tilde{Q}_n$ being distributed as $Q_n$, (b) $\tilde{Q}$ being distributed as $Q$ and (c) $\tilde{Q}_n$ converging to $\tilde{Q}$ almost surely (with respect to $\tilde{P}$) under the topology of uniform convergence on compact sets. Thus, (i) $\tilde{t}_n$, the maximizer of $\tilde{Q}_n$, has the same distribution as $\hat{t}_n$, (ii) $\tilde{t}$, the maximizer of $\tilde{Q}(t)$, has the same distribution as $\arg\max Q(t)$ and (iii) $\tilde{Q}_n(\tilde{t}_n)$ and $\tilde{Q}(\tilde{t})$ have the same distribution as $Q_n(\hat{t}_n)$ and $\max_t Q(t)$, respectively. So it suffices to show that $\tilde{t}_n$ converges in $\tilde{P}^\star$ (outer) probability to $\tilde{t}$ and $\tilde{Q}_n(\tilde{t}_n)$ converges in $\tilde{P}^\star$ (outer) probability to $\tilde{Q}(\tilde{t})$. The convergence of $\tilde{t}_n$ to $\tilde{t}$ in outer probability is shown in Theorem 2.7 of [11].

To show that $\tilde{Q}_n(\tilde{t}_n)$ converges in probability to $\tilde{Q}(\tilde{t})$, we need to show that for fixed $\varepsilon > 0, \delta > 0$, we eventually have

$$P^\star(|\tilde{Q}_n(\tilde{t}_n) - \tilde{Q}(\tilde{t})| > \delta) < \varepsilon.$$

Since $\tilde{t}_n$ and $\tilde{t}$ are $O_p(1)$, there exists $M_\varepsilon > 0$ such that, with

$$A_n^c \equiv \{\tilde{t}_n \notin [-M_\varepsilon, M_\varepsilon]\}, \qquad B_n^c \equiv \{\tilde{t} \notin [-M_\varepsilon, M_\varepsilon]\},$$

$P^\star(A_n^c) < \varepsilon/4$ and $P^\star(B_n^c) < \varepsilon/4$, eventually. Furthermore, as $\tilde{Q}_n$ converges to $\tilde{Q}$ almost surely and therefore in probability, uniformly on every compact set, with

$$C_n^c \equiv \left\{ \sup_{s \in [-M_\varepsilon, M_\varepsilon]} |\tilde{Q}_n(s) - \tilde{Q}(s)| > \delta \right\}$$

we have $P^\star(C_n^c) < \varepsilon/2$, eventually. Hence, $P^\star(A_n^c \cup B_n^c \cup C_n^c) < \varepsilon$, so that $P_\star(A_n \cap B_n \cap C_n) > 1 - \varepsilon$, eventually. But

(5.1) $$A_n \cap B_n \cap C_n \subset \{|\tilde{Q}_n(\tilde{t}_n) - \tilde{Q}(\tilde{t})| \leq \delta\},$$

and consequently

$$P_\star(|\tilde{Q}_n(\tilde{t}_n) - \tilde{Q}(\tilde{t})| \leq \delta) \geq P_\star(A_n \cap B_n \cap C_n) > 1 - \varepsilon,$$

eventually. This implies immediately that

$$P^\star(|\tilde{Q}_n(\tilde{t}_n) - \tilde{Q}(\tilde{t})| > \delta) < \varepsilon$$

for all sufficiently large $n$. It remains to show (5.1). To see this, note that for any $\omega \in A_n \cap B_n \cap C_n$ and $s \in [-M_\varepsilon, M_\varepsilon]$,

$$\tilde{Q}_n(s) = \tilde{Q}(s) + \tilde{Q}_n(s) - \tilde{Q}(s) \leq \tilde{Q}(\tilde{t}) + |\tilde{Q}_n(s) - \tilde{Q}(s)|.$$

Taking the supremum over $s \in [-M_\varepsilon, M_\varepsilon]$ and noting that $\tilde{t}_n \in [-M_\varepsilon, M_\varepsilon]$ on the set $A_n \cap B_n \cap C_n$, we have

$$\tilde{Q}_n(\tilde{t}_n) \leq \tilde{Q}(\tilde{t}) + \sup_{s \in [-M_\varepsilon, M_\varepsilon]} |\tilde{Q}_n(s) - \tilde{Q}(s)|,$$

or equivalently

$$\tilde{Q}_n(\tilde{t}_n) - \tilde{Q}(\tilde{t}) \leq \sup_{s \in [-M_\varepsilon, M_\varepsilon]} |\tilde{Q}_n(s) - \tilde{Q}(s)|.$$

An analogous derivation (replacing $\tilde{Q}_n$ everywhere by $\tilde{Q}$ and $\tilde{t}_n$ by $\tilde{t}$, and vice versa) yields

$$\tilde{Q}(\tilde{t}) - \tilde{Q}_n(\tilde{t}_n) \leq \sup_{s \in [-M_\varepsilon, M_\varepsilon]} |\tilde{Q}(s) - \tilde{Q}_n(s)|.$$

Thus
$$|\tilde{Q}_n(\tilde{t}_n) - \tilde{Q}(\tilde{t})| \leq \sup_{s \in [-M_\varepsilon, M_\varepsilon]} |\tilde{Q}_n(s) - \tilde{Q}(s)| \leq \delta,$$
which completes the proof. $\square$

The following modification of a rate theorem of van der Vaart and Wellner ([21], Theorem 3.2.5) is needed in the proof of Lemma 5.1. The notation $\lesssim$ means that the left-hand side is bounded by a generic constant times the right-hand side.

THEOREM 5.2. *Let $\Theta$ and $\mathcal{F}$ be semimetric spaces. Let $\mathbb{M}_n(\theta, F)$ be stochastic processes indexed by $\theta \in \Theta$ and $F \in \mathcal{F}$. Let $\mathbb{M}(\theta, F)$ be a deterministic function, and let $(\theta_0, F_0)$ be a fixed point in the interior of $\Theta \times \mathcal{F}$. Assume that for every $\theta$ in a neighborhood of $\theta_0$,*

(5.2) $$\mathbb{M}(\theta, F_0) - \mathbb{M}(\theta_0, F_0) \lesssim -d^2(\theta, \theta_0),$$

*where $d(\cdot, \cdot)$ is the semimetric for $\Theta$. Let $\hat{\theta}_n$ be a point of maximum of $\mathbb{M}_n(\theta, \hat{F}_n)$, where $\hat{F}_n$ is random. For each $\varepsilon > 0$, suppose that the following hold:*

(a) *There exists a sequence $\mathcal{F}_{n,\varepsilon}$, $n = 1, 2, \ldots$, of metric subspaces of $\mathcal{F}$, each containing $F_0$ in its interior.*

(b) *For all sufficiently small $\delta > 0$ (say $\delta < \delta_0$, where $\delta_0$ does not depend on $\varepsilon$), and for all sufficiently large $n$,*

(5.3)
$$E^* \sup_{\substack{d(\theta,\theta_0)<\delta \\ F \in \mathcal{F}_{n,\varepsilon}}} |(\mathbb{M}_n(\theta, F) - \mathbb{M}(\theta, F_0)) - (\mathbb{M}_n(\theta_0, F) - \mathbb{M}(\theta_0, F_0))|$$
$$\leq C_\varepsilon \frac{\phi_n(\delta)}{\sqrt{n}}$$

*for a constant $C_\varepsilon > 0$ and functions $\phi_n$ (not depending on $\varepsilon$) such that $\delta \mapsto \phi_n(\delta)/\delta^\alpha$ is decreasing in $\delta$ for some constant $\alpha < 2$ not depending on $n$.*

(c) $P(\hat{F}_n \notin \mathcal{F}_{n,\varepsilon}) < \varepsilon$ *for $n$ sufficiently large.*

*If $r_n^2 \phi_n(r_n^{-1}) \lesssim \sqrt{n}$ for every $n$ and $\hat{\theta}_n \to_p \theta_0$, then $r_n d(\hat{\theta}_n, \theta_0) = O_p(1)$.*

LEMMA 5.1. *The process $Q_n(t)$ defined in the proof of Theorem 2.3 converges in distribution in the space $B_{\text{loc}}(\mathbb{R})$ to the Gaussian process $Q_0(t) \equiv aW(t) - b_0 t^2$, whose distribution is a tight Borel measure concentrated on $C_{\max}(\mathbb{R})$. Here $a$ and $b_0$ are defined in Theorem 2.1. Furthermore, the sequence $\{\hat{t}_n\}$ of maximizers of $\{Q_n(t)\}$ is $O_p(1)$ [and hence converges to $\arg\max_t Q_0(t)$ by Theorem 5.1].*

PROOF. We apply the general approach outlined on page 288 of [21]. Define

$$\mathbb{M}_n(d) = \mathbb{P}_n[(Y - \hat{f}(d^0))(1(X \leq d) - 1(X \leq d^0))],$$
$$\mathbb{M}(d) = P[(Y - f(d^0))(1(X \leq d) - 1(X \leq d^0))].$$

Now, $\tilde{d}_n = \arg\max_{d \in \mathbb{R}} \mathbb{M}_n(d)$ and $d^0 = \arg\max_{d \in \mathbb{R}} \mathbb{M}(d)$ and, in fact, $d^0$ is the unique maximizer of $\mathbb{M}$ under the stipulated conditions. The last assertion needs proof, which will be supplied later. We establish the consistency of $\tilde{d}_n$ for $d^0$ and then find the rate of convergence $r_n$ of $\tilde{d}_n$, in other words, that $r_n$ for which $r_n(\tilde{d}_n - d^0)$ is $O_p(1)$. To establish the consistency of $\tilde{d}_n$ for $d^0$, we apply Corollary 3.2.3 (i) of [21]. We first show that $\sup_{d \in \mathbb{R}} |\mathbb{M}_n(d) - M(d)| \to_p 0$. We can write

$$\sup_{d \in \mathbb{R}} |\mathbb{M}_n(d) - \mathbb{M}(d)|$$
$$\leq \sup_{d \in \mathbb{R}} |(\mathbb{P}_n - P)[(Y - f(d^0))(1(X \leq d) - 1(X \leq d^0))]|$$
$$+ \sup_{d \in \mathbb{R}} |\mathbb{P}_n[(f(d^0) - \hat{f}(d^0))(1(X \leq d) - 1(X \leq d^0))]|.$$

The class of functions $\{(Y - f(d^0))(1(X \leq d) - 1(X \leq d^0)) : d \in \mathbb{R}\}$ is VC with a square integrable envelope [since $E(Y^2) < \infty$] and consequently Glivenko–Cantelli in probability. Thus the first term converges to zero in probability. The second term is easily seen to be bounded by $2|\hat{f}(d^0) - f(d^0)|$, which converges to zero almost surely. It follows that $\sup_{d \in \mathbb{R}} |\mathbb{M}_n(d) - \mathbb{M}(d)| = o_p(1)$. It remains to show that $\mathbb{M}(d^0) > \sup_{d \notin G} \mathbb{M}(d)$ for every open interval $G$ that contains $d^0$. Since $d^0$ is the unique maximizer of the continuous (in fact, differentiable) function $\mathbb{M}(d)$ and $\mathbb{M}(d^0) = 0$, it suffices to show that $\lim_{d \to -\infty} \mathbb{M}(d) < 0$ and $\lim_{d \to \infty} \mathbb{M}(d) < 0$. This is indeed the case, and will be demonstrated at the end of the proof. Thus, all conditions of Corollary 3.2.3 are satisfied, and hence $\tilde{d}_n$ converges in probability to $d^0$.

Next we apply Theorem 5.2 to find the rate of convergence $r_n$ of $\tilde{d}_n$. Given $\varepsilon > 0$, let $\mathcal{F}_{n,\varepsilon} = [f(d^0) - M_\varepsilon/\sqrt{n}, f(d^0) + M_\varepsilon/\sqrt{n}]$, where $M_\varepsilon$ is chosen in such a way that $\sqrt{n}(\hat{f}(d^0) - f(d^0)) \leq M_\varepsilon$, for sufficiently large $n$, with probability at least $1 - \varepsilon$. Since $\hat{f}(d^0) = (\hat{\beta}_l^0 + \hat{\beta}_u^0)/2$ is $\sqrt{n}$-consistent for $f(d^0)$, this can indeed be arranged. Then, setting $\hat{F}_n = \hat{f}(d^0)$, we have $P(\hat{F}_n \notin \mathcal{F}_{n,\varepsilon}) < \varepsilon$ for all sufficiently large $n$. We let $d$ play the role of $\theta$, with $d^0 = \theta_0$, and define

$$\mathbb{M}_n(d, F) = \mathbb{P}_n[(Y - F)(1(X \leq d) - 1(X \leq d^0))],$$
$$\mathbb{M}(d, F) = P[(Y - F)(1(X \leq d) - 1(X \leq d^0))].$$

Then $\tilde{d}_n$ maximizes $\mathbb{M}_n(d, \hat{F}_n) \equiv \mathbb{M}_n(d)$ and $d^0$ maximizes $\mathbb{M}(d, F_0)$, where $F_0 = f(d^0)$. Consequently,

$$\mathbb{M}(d, F_0) - \mathbb{M}(d_0, F_0) \equiv \mathbb{M}(d) - \mathbb{M}(d^0) \leq -C(d - d^0)^2$$

(for some positive constant $C$) for all $d$ in a neighborhood of $d^0$ (say $d \in [d^0 - \delta_0, d^0 + \delta_0]$), on using the continuity of $\mathbb{M}''(d)$ in a neighborhood of $d^0$ and the fact that $\mathbb{M}''(d^0) < 0$ (which follows from arguments at the end of this proof). Thus (5.2) is satisfied. We will next show that (5.3) is also satisfied in our case, with $\phi_n(\delta) \equiv \sqrt{\delta}$, for all $\delta < \delta_0$. Solving $r_n^2 \phi_n(r_n^{-1}) \lesssim \sqrt{n}$ yields $r_n = n^{1/3}$, and we conclude that $n^{1/3}(\tilde{d}_n - d^0) = O_p(1)$.

To show (5.3), we need to find functions $\phi_n(\delta)$ such that

$$E^\star \sup_{|d-d^0|<\delta, F\in\mathcal{F}_{n,\varepsilon}} \sqrt{n}|\mathbb{M}_n(d, F) - \mathbb{M}(d, F_0)|$$

is bounded by $\phi_n(\delta)$. Writing $\mathbb{G}_n \equiv \sqrt{n}(\mathbb{P}_n - P)$, we find that the left-hand side of the above display is bounded by $A_n + B_n$, where

$$A_n = E^\star \sup_{|d-d^0|<\delta, F\in\mathcal{F}_{n,\varepsilon}} |\mathbb{G}_n[(Y-F)(1(X \leq d) - 1(X \leq d^0))]|$$

and

$$B_n = E^\star \sup_{|d-d^0|<\delta, F\in\mathcal{F}_{n,\varepsilon}} \sqrt{n}|P[(F-F_0)(1(X \leq d) - 1(X \leq d^0))]|.$$

First consider the term $A_n$. For sufficiently large $n$,

$$A_n \leq E^\star \sup_{|d-d^0|<\delta, F\in[F_0-1,F_0+1]} |\mathbb{G}_n[(Y-F)(1(X \leq d) - 1(X \leq d^0))]|.$$

Denote by $\mathcal{M}_\delta$ the class of functions $\{(Y-F)(1(X \leq d) - 1(X \leq d^0)): |d - d^0| \leq \delta, F \in [F_0 - 1, F_0 + 1]\}$. An envelope function for this class is given by $M_\delta = (|Y| + F_0 + 2)1(X \in [d^0 - \delta, d^0 + \delta])$. From [21], page 291, using their notation,

$$E^\star(\|\mathbb{G}_n\|_{\mathcal{M}_\delta}) \lesssim J(1, \mathcal{M}_\delta)(PM_\delta^2)^{1/2},$$

where $M_\delta$ is an envelope function for $\mathcal{M}_\delta$ and $J(1, \mathcal{M}_\delta)$ is the uniform entropy integral (considered below). By straightforward computation, there exists $\delta_0 > 0$ such that for all $\delta < \delta_0$, we have $E(M_\delta^2) \lesssim \delta$, for a constant not depending on $\delta$ (but possibly on $\delta_0$). Also, as will be shown below, $J(1, \mathcal{M}_\delta)$ is bounded for all sufficiently small $\delta$. Hence, $A_n \lesssim \sqrt{\delta}$. Next, note that

$$B_n = \sup_{|d-d^0|<\delta, F\in\mathcal{F}_{n,\varepsilon}} \sqrt{n}|P[(F-F_0)(1(X \leq d) - 1(X \leq d^0))]|$$

$$\leq M_\varepsilon \sup_{|d-d^0|<\delta} |F_X(d) - F_X(d^0)| \lesssim M_\varepsilon \delta,$$

CONFIDENCE SETS FOR SPLIT POINTS 23

using condition (A3) in the last step. Hence $A_n + B_n \lesssim \sqrt{\delta} + \delta \lesssim \sqrt{\delta}$, since $\delta$ can be taken less than 1. Thus the choice $\phi_n(\delta) = \sqrt{\delta}$ does indeed work.

Now we check the boundedness of

$$J(1, \mathcal{M}_\delta) = \sup_Q \int_0^1 \sqrt{1 + \log N(\eta \|M_\delta\|_{Q,2}, \mathcal{M}_\delta, L_2(Q))}\, d\eta$$

for small $\delta$, as claimed above. Take any $\eta > 0$. Construct a grid of points on $[F_0 - 1, F_0 + 1]$ such that two successive points on the grid are at distance less than $\eta$ apart. This can be done using fewer than $3/\eta$ points. Now, take a function in $\mathcal{M}_\delta$. This looks like $(Y - F)(1(X \leq d) - 1(X \leq d^0))$ for some $F \in [F_0 - 1, F_0 + 1]$ and some $d$ with $|d - d^0| \leq \delta$. Find the closest point to $F$ on this grid; call this $F_c$. Note that

$$|(Y - F)(1(X \leq d) - 1(X \leq d^0)) - (Y - F_c)(1(X \leq d) - 1(X \leq d^0))|$$
$$\leq \eta 1[X \in [d^0 - \delta, d^0 + \delta]] \leq \eta M_\delta,$$

whence

$$\|(Y - F)(1(X \leq d) - 1(X \leq d^0)) - (Y - F_c)(1(X \leq d) - 1(X \leq d^0))\|_{Q,2}$$

is bounded by $\eta \|M_\delta\|_{Q,2}$. Now for any fixed point $F_{\text{grid}}$ on the grid, $\mathcal{M}_{\delta, F_{\text{grid}}} = \{(Y - F_{\text{grid}})(1(X \leq d) - 1(X \leq d^0)) : d \in [d^0 - \delta, d^0 + \delta]\}$ is a VC-class with VC-dimension bounded by a constant not depending on $\delta$ or $F_{\text{grid}}$. Also, $M_\delta$ is an envelope for $\mathcal{M}_{\delta, F_{\text{grid}}}$; it follows from bounds on covering numbers for VC-classes that $N(\eta \|M_\delta\|_{Q,2}, \mathcal{M}_{\delta, F_{\text{grid}}}, L_2(Q)) \lesssim \eta^{-V_1}$ for some $V_1 > 0$ that does not depend on $Q, F_{\text{grid}}$ or $\delta$. Since the number of grid points is of order $1/\eta$, using the bound on the above display we have

$$N(2\eta \|M_\delta\|_{Q,2}, \mathcal{M}_\delta, L_2(Q)) \lesssim \eta^{-(V_1+1)}.$$

Using this upper bound on the covering number, we obtain a finite upper bound on $J(1, \mathcal{M}_\delta)$ for all $\delta < \delta_0$, via direct computation. This completes the proof that $\hat{t}_n = n^{1/3}(\tilde{d}_n - d^0) = O_p(1)$.

Recalling notation from the proof of Theorem 2.3, we can write

$$Q_n(t) = \xi_n(d^0 + tn^{-1/3}) = R_n(t) + r_{n,1}(t) + r_{n,2}(t),$$

where $R_n(t) = n^{2/3}\, \mathbb{P}_n\, [g(\cdot, d^0 + tn^{-1/3})]$ with

$$g((X,Y), d) = \left(Y - \frac{\beta_l^0 + \beta_u^0}{2}\right)[1(X \leq d) - 1(X \leq d^0)],$$

$$r_{n,1}(t) = n^{1/6}(\hat{f}(d^0) - f(d^0))\sqrt{n}(\mathbb{P}_n - P)[1(X \leq d^0 + tn^{-1/3}) - 1(X \leq d^0)]$$

and

$$r_{n,2}(t) = n^{2/3}(\hat{f}(d^0) - f(d^0))P(1(X \leq d^0 + tn^{-1/3}) - 1(X \leq d^0)).$$

Here, $r_{n,1}(t) \to_p 0$ uniformly on every compact set of the form $[-K, K]$ by applying Donsker's theorem to the empirical process

$$\{\sqrt{n}(\mathbb{P}_n - P)(1(X \leq d^0 + s) - 1(X \leq d^0)) : s \in (-\infty, \infty)\}$$

along with $n^{1/6}(\hat{f}(d^0) - f(d^0)) = o_p(1)$. The term $r_{n,2}(t) \to_p 0$ uniformly on every $[-K, K]$ since $n^{1/3}(\hat{f}(d^0) - f(d^0)) = o_p(1)$ and $n^{1/3} \sup_{t \in [-K,K]} P(1(X \leq d^0 + tn^{-1/3}) - 1(X \leq d^0)) = O(1)$. Hence, the limiting distribution of $Q_n(t)$ will be the same as the limiting distribution of $R_n(t)$. We show that $R_n \to_d Q_0$, where $Q_0$ is the Gaussian process defined in Theorem 2.3. Write

$$R_n(t) = n^{2/3}(\mathbb{P}_n - P)[g(\cdot, d^0 + tn^{-1/3})] + n^{2/3} P[g(\cdot, d^0 + tn^{-1/3})]$$
$$= I_n(t) + J_n(t).$$

In terms of the empirical process $\mathbb{G}_n$, we have $I_n(t) = \mathbb{G}_n(f_{n,t})$ where

$$f_{n,t}(x, y) = n^{1/6}(y - f(d^0))(1(x \leq d^0 + tn^{-1/3}) - 1(x \leq d^0)).$$

We will use Theorem 2.11.22 from [21] to show that on each compact set $[-K, K]$, $\mathbb{G}_n f_{n,t}$ converges as a process in $l^\infty[-K, K]$ to the tight Gaussian process $aW(t)$, where $a^2 = \sigma^2(d^0) p_X(d^0)$. Also, $J_n(t)$ converges on every $[-K, K]$ uniformly to the deterministic function $-b_0 t^2$, with $b_0 = |f'(d^0)| p_X(d^0)/2 > 0$. Hence $Q_n(t) \to_d Q_0(t) \equiv aW(t) - b_0 t^2$ in $B_{\text{loc}}(\mathbb{R})$, as required.

To complete the proof, we need to show that $I_n$ and $J_n$ have the limits claimed above. As far as $I_n$ is concerned, provided we can verify the other conditions of Theorem 2.11.22 from [21], the covariance kernel $H(s,t)$ of the limit of $\mathbb{G}_n f_{n,t}$ is given by the limit of $P(f_{n,s} f_{n,t}) - P f_{n,s} P f_{n,t}$ as $n \to \infty$. We first compute $P(f_{n,s} f_{n,t})$. This vanishes if $s$ and $t$ are of opposite signs. For $s, t > 0$,

$$Pf_{n,s} f_{n,t} = E[n^{1/3}(Y - f(d^0))^2 1\{X \in (d^0, d^0 + (s \wedge t) n^{-1/3}]\}]$$
$$= \int_{d^0}^{d^0 + (s \wedge t) n^{-1/3}} n^{1/3} [E[(f(X) + \varepsilon - f(d^0))^2 | X = x]] p_X(x) \, dx$$
$$= n^{1/3} \int_{d^0}^{d^0 + (s \wedge t) n^{-1/3}} (\sigma^2(x) + (f(x) - f(d^0))^2) p_X(x) \, dx$$
$$\to \sigma^2(d^0) p_X(d^0)(s \wedge t)$$
$$\equiv a^2(s \wedge t).$$

Also, it is easy to see that $Pf_{n,s}$ and $Pf_{n,t}$ converge to 0. Thus, when $s, t > 0$,

$$P(f_{n,s} f_{n,t}) - Pf_{n,s} Pf_{n,t} \to a^2(s \wedge t) \equiv H(s,t).$$

Similarly, it can be checked that for $s, t < 0$, $H(s,t) = a^2(-s \wedge -t)$. Thus $H(s,t)$ is the covariance kernel of the Gaussian process $aW(t)$.

Next we need to check

$$\sup_Q \int_0^{\delta_n} \sqrt{\log N(\varepsilon\|F_n\|_{Q,2}, \mathcal{F}_n, L_2(Q))}\, d\varepsilon \to 0, \tag{5.4}$$

for every $\delta_n \to 0$, where

$$\mathcal{F}_n = \{n^{1/6}(y - f(d^0))[1(x \le d^0 + tn^{-1/3}) - 1(x \le d^0)] : t \in [-K, K]\}$$

and

$$F_n(x, y) = n^{1/6}|y - f(d^0)|1(x \in [d^0 - Kn^{-1/3}, d^0 + Kn^{-1/3}])$$

is an envelope for $\mathcal{F}_n$. From [21], page 141,

$$N(\varepsilon\|F_n\|_{Q,2}, \mathcal{F}_n, L_2(Q)) \le KV(\mathcal{F}_n)(16e)^{V(\mathcal{F}_n)}\left(\frac{1}{\varepsilon}\right)^{2(V(\mathcal{F}_n)-1)}$$

for a universal constant $K$ and $0 < \varepsilon < 1$, where $V(\mathcal{F}_n)$ is the VC-dimension of $\mathcal{F}_n$. Since $V(\mathcal{F}_n)$ is uniformly bounded, we see that the above inequality implies

$$N(\varepsilon\|F_n\|_{Q,2}, \mathcal{F}_n, L_2(Q)) \lesssim \left(\frac{1}{\varepsilon}\right)^s,$$

where $s = \sup_n 2(V(\mathcal{F}_n) - 1) < \infty$, so (5.4) follows from

$$\int_0^{\delta_n} \sqrt{-\log \varepsilon}\, d\varepsilon \to 0$$

as $\delta_n \to 0$. We also need to check the conditions (2.11.21) in [21]:

$$P^\star F_n^2 = O(1), \qquad P^\star F_n^2 1\{F_n > \eta\sqrt{n}\} \to 0 \qquad \forall \eta > 0,$$

and

$$\sup_{|s-t|<\delta_n} P(f_{n,s} - f_{n,t})^2 \to 0 \qquad \forall \delta_n \to 0.$$

With $F_n$ as defined above, an easy computation shows that

$$P^\star F_n^2 = K\frac{1}{Kn^{-1/3}}\int_{d^0-Kn^{-1/3}}^{d^0+Kn^{-1/3}}(\sigma^2(x) + (f(x)-f(d^0))^2)p_X(x)\,dx = O(1).$$

Denote the set $[d^0 - Kn^{-1/3}, d^0 + Kn^{-1/3}]$ by $S_n$. Then

$$\begin{aligned}
&P^\star(F_n^2 1\{F_n > \eta\sqrt{n}\}) \\
&\quad = E[n^{1/3}|Y-f(d^0)|^2 1\{X \in S_n\}1\{|Y-f(d^0)|1\{X \in S_n\} > \eta n^{1/3}\}] \\
&\quad \le E[n^{1/3}|Y-f(d^0)|^2 1\{X \in S_n\}1\{|\varepsilon| > \eta n^{1/3}/2\}] \\
&\quad \le E[2n^{1/3}(\varepsilon^2 + (f(X)-f(d^0))^2)1\{X \in S_n\}1\{|\varepsilon| > \eta n^{1/3}/2\}],
\end{aligned} \tag{5.5}$$

eventually, since for all sufficiently large $n$

$$\{|Y - f(d^0)|1\{X \in S_n\} > \eta n^{1/3}\} \subset \{|\varepsilon| > \eta n^{1/3}/2\}.$$

Now, the right-hand side of (5.5) can be written as $T_1 + T_2$, where

$$T_1 = 2n^{1/3} E[\varepsilon^2 1\{|\varepsilon| > \eta n^{1/3}/2\}1\{X \in S_n\}]$$

and

$$T_2 = 2n^{1/3} E[(f(X) - f(d^0))^2 1\{X \in S_n\}1\{|\varepsilon| > \eta n^{1/3}/2\}].$$

We will show that $T_1 = o(1)$. We have

$$T_1 = 2n^{1/3} \int_{d^0 - Kn^{-1/3}}^{d^0 + Kn^{-1/3}} E[\varepsilon^2 1\{|\varepsilon| > \eta n^{1/3}/2\}|X = x] p_X(x)\, dx.$$

By (A5), for any $\xi > 0$,

$$\sup_{x \in S_n} E[\varepsilon^2 1\{|\varepsilon| > \eta n^{1/3}/2\}|X = x] < \xi$$

for $n$ sufficiently large. Since $n^{1/3} \int_{S_n} p_X(x)\, dx$ is eventually bounded by $2 \times K p_X(d^0)$ it follows that $T_1$ is eventually smaller than $2\xi K p_X(d^0)$. We conclude that $T_1 = o(1)$. Next, note that (A5) implies that $\sup_{x \in S_n} E[1\{|\varepsilon| > \eta n^{1/3}/2\}|X = x] \to 0$ as $\eta \to \infty$, so $T_2 = o(1)$ by an argument similar to that above. Finally,

$$\sup_{|s-t|<\delta_n} P(f_{n,s} - f_{n,t})^2 \to 0$$

as $\delta_n \to 0$ can be checked via similar computations.

We next deal with $J_n$. For convenience we sketch the uniformity of the convergence of $J_n(t)$ to the claimed limit on $0 \leq t \leq K$. We have

$$\begin{aligned}
J_n(t) &= n^{2/3} E[(Y - f(d^0))(1(X \leq d^0 + tn^{-1/3}) - 1(X \leq d^0))] \\
&= n^{2/3} E[(f(X) - f(d^0))1(X \in (d^0, d^0 + tn^{-1/3}])] \\
&= n^{2/3} \int_{d^0}^{d^0 + tn^{-1/3}} (f(x) - f(d^0)) p_X(x)\, dx \\
&= n^{1/3} \int_0^t (f(d^0 + un^{-1/3}) - f(d^0)) p_X(d^0 + un^{-1/3})\, du \\
&= \int_0^t u \frac{f(d^0 + un^{-1/3}) - f(d^0)}{un^{-1/3}} p_X(d^0 + un^{-1/3})\, du \\
&\to \int_0^t u f'(d^0) p_X(d^0)\, du \qquad \text{(uniformly on } 0 \leq t \leq K\text{)} \\
&= \frac{1}{2} f'(d^0) p_X(d^0) t^2.
\end{aligned}$$

It only remains to verify that (i) $d^0$ is the unique maximizer of $\mathbb{M}(d)$, (ii) $\mathbb{M}(-\infty) < 0, \mathbb{M}(\infty) < 0$ and (iii) $f'(d^0)p_X(d^0) < 0$ [so the process $aW(t) + (f'(d^0)p_X(d^0)/2)t^2$ is indeed in $C_{\max}(\mathbb{R})$]. To show (i), recall that

$$\mathbb{M}(d) = E[g((X,Y),d)] = E\left[\left(Y - \frac{\beta_l^0 + \beta_u^0}{2}\right)(1(X \leq d) - 1(X \leq d^0))\right].$$

Let $\xi(d) = E[Y - \beta_l^0 1(X \leq d) - \beta_u^0 1(X > d)]^2$. By condition (A1), $d^0$ is the unique minimizer of $\xi(d)$. Consequently, $d^0$ is also the unique maximizer of the function $\xi(d^0) - \xi(d)$. Straightforward algebra shows that

$$\xi(d^0) - \xi(d) = 2(\beta_l^0 - \beta_u^0)\mathbb{M}(d)$$

and since $\beta_l^0 - \beta_u^0 > 0$, it follows that $d^0$ is also the unique maximizer of $\mathbb{M}(d)$. This shows (i). Next,

$$\mathbb{M}(d) = E[(f(X) - f(d^0))(1(X \leq d) - 1(X \leq d^0))]$$
$$= \int_{-\infty}^{\infty} (f(x) - f(d^0))(1(x \leq d) - 1(x \leq d^0))p_X(x)\,dx$$
$$= \int_{-\infty}^{d} (f(x) - f(d^0))p_X(x)\,dx - \int_{-\infty}^{d^0} (f(x) - f(d^0))p_X(x)\,dx,$$

so that

$$\mathbb{M}(-\infty) = \lim_{d \to -\infty} \mathbb{M}(d) = -\int_{-\infty}^{d^0} (f(x) - f(d^0))p_X(x)\,dx < 0$$

if and only if $\int_{-\infty}^{d^0} f(x)p_X(x)\,dx > f(d^0)F_X(d^0)$ if and only if $\beta_l^0 \equiv \int_{-\infty}^{d^0} f(x) \times p_X(x)\,dx/F_X(d^0) > (\beta_l^0 + \beta_u^0)/2$, and this is indeed the case, since $\beta_l^0 > \beta_u^0$. We can prove that $\mathbb{M}(\infty) < 0$ in a similar way, so (ii) holds. Also, $\mathbb{M}'(d) = (f(d) - f(d^0))p_X(d)$, so $\mathbb{M}'(d^0) = 0$. Finally,

$$\mathbb{M}''(d) = f'(d)p_X(d) + (f(d) - f(d^0))p_X'(d),$$

so $\mathbb{M}''(d^0) = f'(d^0)p_X(d^0) \leq 0$, since $d^0$ is the maximizer. This implies (iii), since by our assumptions $f'(d^0)p_X(d^0) \neq 0$. □

PROOF OF THEOREM 2.1. Let $\Theta$ denote the set of all possible values of $(\beta_l, \beta_u, d)$ and let $\theta$ denote a generic vector in $\Theta$. Define the criterion function $\mathbb{M}(\theta) = Pm_\theta$, where

$$m_\theta(x,y) = (y - \beta_l)^2 1(x \leq d) + (y - \beta_u)^2 1(x > d).$$

The vector $\theta_0 \equiv (\beta_l^0, \beta_u^0, d^0)$ minimizes $\mathbb{M}(\theta)$, while $\hat{\theta}_n \equiv (\hat{\beta}_l, \hat{\beta}_u, \hat{d}_n)$ minimizes $\mathbb{M}_n(\theta) = \mathbb{P}_n m_\theta$. Since $\theta_0$ uniquely minimizes $\mathbb{M}(\theta)$ under condition (A1), using the twice continuous differentiability of $\mathbb{M}$ at $\theta_0$, we have

$$\mathbb{M}(\theta) - \mathbb{M}(\theta_0) \geq \tilde{C}d^2(\theta, \theta_0)$$

in a neighborhood of $\theta_0$ (for some $\tilde{C} > 0$), where $d(\cdot, \cdot)$ is the $l_\infty$ metric on $\mathbb{R}^3$. Thus, there exists $\delta_0 > 0$ sufficiently small, such that for all $(\beta_l, \beta_u, d)$ with $|\beta_l - \beta_l^0| < \delta_0$, $|\beta_u - \beta_u^0| < \delta_0$ and $|d - d^0| < \delta_0$, the above display holds.

For all $\delta < \delta_0$ we will find a bound on $E_P^\star \|\mathbb{G}_n\|_{\mathcal{M}_\delta}$, where $\mathcal{M}_\delta \equiv \{m_\theta - m_{\theta_0} : d(\theta, \theta_0) < \delta\}$ and $\mathbb{G}_n \equiv \sqrt{n}(\mathbb{P}_n - P)$. From [21], page 298,

$$E_P^\star \|\mathbb{G}_n\|_{\mathcal{M}_\delta} \leq J(1, \mathcal{M}_\delta)(PM_\delta^2)^{1/2},$$

where $M_\delta$ is an envelope function for the class $\mathcal{M}_\delta$. Straightforward algebra shows that

$$\begin{aligned}(m_\theta &- m_{\theta_0})(X, Y) \\ &= 2(Y - f(d^0))(\beta_u^0 - \beta_l^0)\{1(X \leq d) - 1(X \leq d^0)\} \\ &\quad + (\beta_l^0 - \beta_l)(2Y - \beta_l^0 - \beta_l)1(X \leq d) \\ &\quad + (\beta_u^0 - \beta_u)(2Y - \beta_u^0 - \beta_u)1(X > d).\end{aligned}$$

The class of functions

$$\mathcal{M}_{1,\delta} = \{2(Y - f(d^0))(\beta_u^0 - \beta_l^0)\{1(X \leq d) - 1(X \leq d^0)\} : d \in [d^0 - \delta, d^0 + \delta]\}$$

is easily seen to be VC, with VC-dimension bounded by a constant not depending on $\delta$; furthermore, $M_{1,\delta} = 2|(Y - f(d^0))(\beta_u^0 - \beta_l^0)|1(X \in [d^0 - \delta, d^0 + \delta])$ is an envelope function for this class. It follows that

$$N(\varepsilon \|M_{1,\delta}\|_{P,2}, \mathcal{M}_{1,\delta}, L_2(P)) \lesssim \varepsilon^{-V_1},$$

for some $V_1 > 0$ that does not depend on $\delta$. Next, consider the class of functions

$$\begin{aligned}\mathcal{M}_{2,\delta} = \{(\beta_l^0 - \beta_l)(2Y - \beta_l^0 - \beta_l)1(X \leq d) &: d \in [d^0 - \delta, d^0 + \delta], \\ &\beta_l \in [\beta_l^0 - \delta, \beta_l^0 + \delta]\}.\end{aligned}$$

Fix a grid of points $\{\beta_{l,c}\}$ in $[\beta_l^0 - \delta, \beta_l^0 + \delta]$ such that successive points on this grid are at a distance less than $\tilde{\varepsilon}$ apart, where $\tilde{\varepsilon} = \varepsilon\delta/2$. The cardinality of this grid is certainly less than $3\delta/\tilde{\varepsilon}$. For a fixed $\beta_{l,c}$ in this grid, the class of functions $\mathcal{M}_{2,\delta,c} \equiv \{(\beta_l^0 - \beta_{l,c})(2Y - \beta_l^0 - \beta_{l,c})1(X \leq d) : d \in [d^0 - \delta, d^0 + \delta]\}$ is certainly VC with VC-dimension bounded by a constant that does not depend on $\delta$ or the point $\beta_{l,c}$. Also, note that $M_{2,\delta} \equiv \delta(2|Y| + C)$, where $C$ is a sufficiently large constant not depending on $\delta$, is an envelope function for the class $\mathcal{M}_{2,\delta}$, and hence also an envelope function for the restricted class with $\beta_{l,c}$ held fixed. It follows that for some universal positive constant $V_2 > 0$ and any $\eta > 0$,

$$N(\eta \|M_{2,\delta}\|_{P,2}, \mathcal{M}_{2,\delta,c}, L_2(P)) \lesssim \eta^{-V_2}.$$

Now, $\|M_{2,\delta}\|_{P,2} = \delta\|G\|_{P,2}$, where $G = 2|Y| + C$. Thus,

$$N(\tilde{\varepsilon}\|G\|_{P,2}, \mathcal{M}_{2,\delta,c}, L_2(P)) \lesssim \left(\frac{\delta}{\tilde{\varepsilon}}\right)^{V_2}.$$

Next, consider a function $g(X,Y) = (\beta_l^0 - \beta_l)(2Y - \beta_l^0 - \beta_l)1(X \leq d)$ in $\mathcal{M}_{2,\delta}$. Find a $\beta_{l,c}$ that is within $\tilde{\varepsilon}$ distance of $\beta_l$. There are of order $(\delta/\tilde{\varepsilon})^{V_2}$ balls of radius $\tilde{\varepsilon}\|G\|_{P,2}$ that cover the class $\mathcal{M}_{2,\delta,c}$, so the function $g_c(X,Y) \equiv (\beta_l^0 - \beta_{l,c})(2Y - \beta_l^0 - \beta_{l,c})1(X \leq d)$ must be at distance less than $\tilde{\varepsilon}\|G\|_{P,2}$ from the center of one of these balls, say $B$. Also, it is easily checked that $\|g - g_c\|_{P,2} < \tilde{\varepsilon}\|G\|_{P,2}$. Hence $g$ must be at distance less than $2\tilde{\varepsilon}\|G\|_{P,2}$ from the center of $B$. It then readily follows that

$$N(2\tilde{\varepsilon}\|G\|_{P,2}, \mathcal{M}_{2,\delta}, L_2(P)) \lesssim \left(\frac{\delta}{\tilde{\varepsilon}}\right)^{V_2+1},$$

on using the fact that the cardinality of the grid $\{\beta_{l,c}\}$ is of order $\delta/\tilde{\varepsilon}$. Substituting $\varepsilon\delta/2$ for $\tilde{\varepsilon}$ in the above display, we get

$$N(\varepsilon\|M_{2,\delta}\|_{P,2}, \mathcal{M}_{2,\delta}, L_2(P)) \lesssim \left(\frac{1}{\varepsilon}\right)^{V_2+1}.$$

Finally, with

$$\mathcal{M}_{3,\delta} = \{(\beta_u^0 - \beta_u)(2Y - \beta_u^0 - \beta_u)1(X > d) : d \in [d^0 - \delta, d^0 + \delta],$$
$$\beta_u \in [\beta_u^0 - \delta, \beta_u^0 + \delta]\}$$

and $M_{3,\delta} = \delta(2|Y| + C')$ for some sufficiently large constant $C'$ not depending on $\delta$, we similarly argue that

$$N(\varepsilon\|M_{3,\delta}\|_{P,2}, \mathcal{M}_{3,\delta}, L_2(P)) \lesssim \left(\frac{1}{\varepsilon}\right)^{V_3+1},$$

for some positive constant $V_3$ not depending on $\delta$. The class $\mathcal{M}_\delta \subset \mathcal{M}_{1,\delta} + \mathcal{M}_{2,\delta} + \mathcal{M}_{3,\delta} \equiv \tilde{\mathcal{M}}_\delta$. Set $M_\delta = M_{1,\delta} + M_{2,\delta} + M_{3,\delta}$. Now, it is not difficult to see that

$$N(3\varepsilon\|M_\delta\|_{P,2}, \tilde{\mathcal{M}}_\delta, L_2(P)) \lesssim \left(\frac{1}{\varepsilon}\right)^{V_1+V_2+V_3}.$$

This also holds for any probability measure $Q$ such that $0 < E_Q(Y^2) < \infty$, with the constant being independent of $Q$ or $\delta$. Since $\mathcal{M}_\delta \subset \tilde{\mathcal{M}}_\delta$, it follows that

$$N(3\varepsilon\|M_\delta\|_{Q,2}, \mathcal{M}_\delta, L_2(Q)) \lesssim \left(\frac{1}{\varepsilon}\right)^{V_1+V_2+V_3}.$$

Thus, with $\mathcal{Q}$ denoting the set of all such measures $Q$,

$$J(1, M_\delta) \equiv \sup_{Q \in \mathcal{Q}} \int_0^1 \sqrt{1 + \log N(\varepsilon\|M_\delta\|_{Q,2}, \mathcal{M}_\delta, L_2(Q))}\, d\varepsilon < \infty$$

for all sufficiently small $\delta$. Next,
$$PM_\delta^2 \lesssim PM_{1,\delta}^2 + PM_{2,\delta}^2 + PM_{3,\delta}^2 \lesssim \delta + \delta^2 \lesssim \delta$$
since we can assume $\delta < 1$. Therefore $E_P^\star \|\mathbb{G}_n\|_{\mathcal{M}_\delta} \lesssim \sqrt{\delta}$, and $\phi_n(\delta)$ in Theorem 3.2.5 of [21] can be taken as $\sqrt{\delta}$. Solving $r_n^2 \phi_n(1/r_n) \leq \sqrt{n}$ yields $r_n \leq n^{1/3}$, and we conclude that
$$n^{1/3}(\hat{\beta}_l - \beta_l^0, \hat{\beta}_u - \beta_u^0, \hat{d}_n - d^0) = O_p(1).$$

Having established the rate of convergence, we now determine the asymptotic distribution. It is easy to see that
$$n^{1/3}(\hat{\beta}_l - \beta_l^0, \hat{\beta}_u - \beta_u^0, \hat{d}_n - d^0) = \arg\min_h V_n(h),$$
where
$$(5.6)\quad V_n(h) = n^{2/3}(\mathbb{P}_n - P)[m_{\theta_0 + hn^{-1/3}} - m_{\theta_0}] + n^{2/3}P[m_{\theta_0 + hn^{-1/3}} - m_{\theta_0}]$$
for $h = (h_1, h_2, h_3) \in \mathbb{R}^3$. The second term above converges to $h^T V h/2$, uniformly on every $[-K, K]^3$ ($K > 0$), where $V$ is the Hessian of the function $\theta \mapsto Pm_\theta$ at the point $\theta_0$, on using the twice continuous differentiability of the function at $\theta_0$ and the fact that $\theta_0$ minimizes this function. Note that $V$ is a positive definite matrix. Calculating the Hessian matrix gives
$$V = \begin{pmatrix} 2F_X(d^0) & 0 & (\beta_l^0 - \beta_u^0)p_X(d^0) \\ 0 & 2(1 - F_X(d^0)) & (\beta_l^0 - \beta_u^0)p_X(d^0) \\ (\beta_l^0 - \beta_u^0)p_X(d^0) & (\beta_l^0 - \beta_u^0)p_X(d^0) & 2|(\beta_l^0 - \beta_u^0)f'(d^0)\,p_X(d^0)| \end{pmatrix}.$$

We next deal with distributional convergence of the first term in (5.6), which can be written as $\sqrt{n}(\mathbb{P}_n - P) f_{n,h}$, where $f_{n,h} = f_{n,h,1} + f_{n,h,2} + f_{n,h,3}$ and
$$f_{n,h,1}(x,y) = n^{1/6} 2(\beta_u^0 - \beta_l^0)(y - f(d^0))(1(x \leq d^0 + h_3 n^{-1/3}) - 1(x \leq d^0)),$$
$$f_{n,h,2}(x,y) = -n^{-1/6} h_1 (2y - 2\beta_l^0 - h_1 n^{-1/3}) 1(x \leq d_0 + h_3 n^{-1/3}),$$
$$f_{n,h,3}(x,y) = -n^{-1/6} h_2 (2y - 2\beta_u^0 - h_2 n^{-1/3}) 1(x > d_0 + h_3 n^{-1/3}).$$

A natural envelope function $F_n$ for $\mathcal{F}_n \equiv \{f_{n,h} : h \in [-K, K]^3\}$ is given by
$$F_n(x, y) = 2n^{1/6}|(\beta_l^0 - \beta_u^0)(y - f(d^0))|1\{x \in [d^0 - Kn^{-1/3}, d^0 + Kn^{-1/3}]\}$$
$$+ Kn^{-1/6}(2|y - \beta_l^0| + 1) + Kn^{-1/6}(2|y - \beta_u^0| + 1).$$

The limiting distribution of $\sqrt{n}(\mathbb{P}_n - P) f_{n,h}$ is directly obtained by appealing to Theorem 2.11.22 of [21]. On each compact set of the form $[-K, K]^3$, the process $\sqrt{n}(\mathbb{P}_n - P) f_{n,h}$ converges in distribution to $\tilde{a} W(h_3)$, where $\tilde{a} = 2|\beta_l^0 - \beta_u^0|(\sigma^2(d^0) p_X(d^0))^{1/2}$. This follows on noting that
$$\lim_{n \to \infty} P f_{n,s} f_{n,h} - P f_{n,s} P f_{n,h} = \tilde{a}^2 (|s_3| \wedge |h_3|) 1(s_3 h_3 > 0),$$

by direct computation and verification of conditions (2.11.21) preceding the statement of Theorem 2.11.22; we omit the details as they are similar to those in the proof of Lemma 5.1. The verification of the entropy-integral condition, that is,

$$\sup_Q \int_0^{\delta_n} \sqrt{\log N(\varepsilon \|F_n\|_{Q,2}, \mathcal{F}_n, L_2(Q))}\, d\varepsilon \to 0$$

as $\delta_n \to 0$, uses $N(\varepsilon \|F_n\|_{Q,2}, \mathcal{F}_n, L_2(Q)) \lesssim \varepsilon^{-V}$ for some $V > 0$ not depending on $Q$; the argument is similar to the one we used earlier with $J(1, \mathcal{M}_\delta)$.

It follows that the process $V_n(h)$ converges in distribution in the space $B_{\mathrm{loc}}(\mathbb{R}^3)$ to the process $\tilde{W}(h_1, h_2, h_3) \equiv \tilde{a} W(h_3) + h^T V h/2$. The limiting distribution is concentrated on $C_{\min}(\mathbb{R}^3)$ [defined analogously to $C_{\max}(\mathbb{R}^3)$], which follows on noting that the covariance kernel of the Gaussian process $\tilde{W}$ has the rescaling property (2.4) of [11] and that $V$ is positive definite; furthermore, $\tilde{W}(s) - \tilde{W}(h)$ has nonzero variance for $s \neq h$, whence Lemma 2.6 of [11] forces a unique minimizer. Invoking Theorem 5.1 (to be precise, a version of the theorem with max replaced by min), we conclude that

$$(5.7) \quad \left(\arg\min_h V_n(h), \min_h V_n(h)\right) \xrightarrow{d} \left(\arg\min_h \tilde{W}(h), \min_h \tilde{W}(h)\right).$$

But note that

$$\min_h \tilde{W}(h) = \min_{h_3}\left\{\tilde{a} W(h_3) + \min_{h_1, h_2} h^T V h/2\right\}$$

and we can find $\arg\min_{h_1, h_2} h^T V h/2$ explicitly. After some routine calculus, we find that the limiting distribution of the first component in (5.7) can be expressed in the form stated in the theorem. This completes the proof. $\square$

PROOF OF THEOREM 2.2. Inspecting the second component of (5.7), we find

$$n^{-1/3} \mathrm{RSS}_0(\beta_l^0, \beta_u^0, d^0) = -\min_h V_n(h) \xrightarrow{d} -\min_h \tilde{W}(h)$$

and this simplifies to the limit stated in the theorem. To show that $n^{-1/3} \mathrm{RSS}_1(d^0)$ converges to the same limit, it suffices to show that the difference $D_n = n^{-1/3} \mathrm{RSS}_0(\beta_l^0, \beta_u^0, d^0) - n^{-1/3} \mathrm{RSS}_1(d^0)$ is asymptotically negligible. Some algebra gives that $D_n = I_n + J_n$, where

$$I_n = n^{-1/3} \sum_{i=1}^n (2Y_i - \hat{\beta}_l^0 - \beta_l^0)(\hat{\beta}_l^0 - \beta_l^0) \mathbf{1}(X_i \leq d^0)$$

and

$$J_n = n^{-1/3} \sum_{i=1}^n (2Y_i - \hat{\beta}_u^0 - \beta_u^0)(\hat{\beta}_u^0 - \beta_u^0) \mathbf{1}(X_i > d^0).$$

Then
$$I_n = \sqrt{n}(\hat{\beta}_l^0 - \beta_l^0)n^{1/6}\mathbb{P}_n[(2Y - \hat{\beta}_l^0 - \beta_l^0)1(X \le d^0)]$$
$$= \sqrt{n}(\hat{\beta}_l^0 - \beta_l^0)n^{1/6}(\mathbb{P}_n - P)[(2Y - \hat{\beta}_l^0 - \beta_l^0)1(X \le d^0)]$$
$$+ \sqrt{n}(\hat{\beta}_l^0 - \beta_l^0)n^{1/6}P[(2Y - \hat{\beta}_l^0 - \beta_l^0)1(X \le d^0)]$$
$$= I_{n,1} + I_{n,2}.$$

Since $\sqrt{n}(\hat{\beta}_l^0 - \beta_l^0) = O_p(1)$ and
$$n^{1/6}(\mathbb{P}_n - P)[(2Y - \hat{\beta}_l^0 - \beta_l^0)1(X \le d^0)]$$
$$= n^{1/6}(\mathbb{P}_n - P)[(2Y - \beta_l^0)1(X \le d^0)]$$
$$- \hat{\beta}_l^0 n^{1/6}(\mathbb{P}_n - P)(1(X \le d^0))$$

is clearly $o_p(1)$ by the central limit theorem and the consistency of $\hat{\beta}_l^0$, we have that $I_{n,1} = o_p(1)$. To show $I_{n,2} = o_p(1)$, it suffices to show that $n^{1/6}P[(2Y - \hat{\beta}_l^0 - \beta_l^0)1(X \le d^0)] \to 0$. But this can be written as
$$n^{1/6}P[2(Y - \beta_l^0)1(X \le d^0)] + n^{1/6}P[(\beta_l^0 - \hat{\beta}_l^0)1(X \le d^0)].$$

The first term vanishes, from the normal equations characterizing $(\beta_l^0, \beta_u^0, d^0)$, and the second term is $n^{1/6}O(n^{-1/2}) \to 0$. We have shown that $I_n = o_p(1)$, and $J_n = o_p(1)$ can be shown in the same way. This completes the proof. $\square$

**Acknowledgments.** The authors thank Song Qian for comments about the Everglades application, Michael Woodroofe and Bin Yu for helpful discussion, Marloes Maathuis for providing the extended rate convergence theorem, and the referees for their detailed comments.

Department of Statistics  
University of Michigan  
1085 South University  
Ann Arbor, Michigan 48109-1107  
USA  
E-mail: moulib@umich.edu

Department of Biostatistics  
Columbia University  
722 West 168th Street, 6th Floor  
New York, New York 10032-2603  
USA  
E-mail: im2131@columbia.edu